\documentclass[11pt]{amsart} \textwidth=14.5cm \oddsidemargin=1cm \evensidemargin=1cm \usepackage{amsmath}
\usepackage{amsxtra} \usepackage{amscd} \usepackage{amsthm} \usepackage{amsfonts} \usepackage{amssymb}
\usepackage{eucal}

\numberwithin{equation}{section} \newtheorem{thm}{Theorem} \newtheorem{cor}[thm]{Corollary}
\newtheorem{lem}[thm]{Lemma} \newtheorem{prop}[thm]{Proposition} 

\theoremstyle{definition} 

\theoremstyle{remark} \newtheorem{rem}[thm]{Remark} \newtheorem*{rem*}{Remark}

\numberwithin{thm}{section}

\theoremstyle{definition}

\theoremstyle{remark}

\newcommand{\nc}{\newcommand} \nc{\renc}{\renewcommand} \nc{\ssec}{\subsection} \nc{\sssec}{\subsubsection}
\nc{\on}{\operatorname} \nc\wt{\widetilde}

\newcommand{\iso}{\stackrel{\sim}{\longrightarrow}} \newcommand{\CY}{{\mathcal Y}}
\newcommand{\CalD}{{\mathcal D}} \newcommand{\Char}{{\mathcal{CS}}} \newcommand{\HC}{{\widehat{HC}}}
\newcommand{\hc}{{\widehat{hc}}} \newcommand{\Har}{{\mathcal{HC}}} \newcommand{\CH}{{\widehat{CH}}}
  
\newcommand{\Lotimes}{{\stackrel{L}{\otimes}}}

\nc{\BA}{{\mathbb{A}}} \nc{\BC}{{\mathbb{C}}} \nc{\BD}{{\mathbb{D}}} \nc{\BF}{{\mathbb{F}}}
\nc{\BM}{{\mathbb{M}}} \nc{\BN}{{\mathbb{N}}} \nc{\BP}{{\mathbb{P}}} \nc{\BR}{{\mathbb{R}}}
\nc{\BZ}{{\mathbb{Z}}} \nc{\BS}{{\mathbb{S}}}

\nc{\Zet}{{\mathbb{Z}}} \nc{\Ce}{{\mathbb{C}}} \nc{\Fq}{{\mathbb{F}_q}}

\nc{\CA}{{\mathcal{A}}} \nc{\A}{{\mathcal{A}}} \nc{\CB}{{\mathcal{B}}} \nc{\B}{{\mathcal{B}}}
\nc{\CC}{{\mathcal{C}}} 
\nc{\CE}{{\mathcal{E}}} \nc{\CF}{{\mathcal{F}}}
\nc{\F}{{\mathcal{F}}} \nc{\CG}{{\mathcal{G}}} \nc{\G}{{\mathcal{G}}} \nc{\CK}{{\mathcal{K}}}
\nc{\CI}{{\mathcal{I}}} \nc{\CJ}{{\mathcal J}} \nc{\CL}{{\mathcal{L}}} \nc{\CM}{{\mathcal{M}}}
\nc{\CMM}{{\mathcal{M}^{\operatorname{gen}}_\hbar(-\rho)}} \nc{\CN}{{\mathcal{N}}} \nc{\CO}{{\mathcal{O}}}
\nc{\CP}{{\mathcal{P}}} \nc{\CQ}{{\mathcal{Q}}} \nc{\CR}{{\mathcal{R}}} \nc{\CS}{{\mathcal{S}}}
\nc{\CT}{{\mathcal{T}}} \nc{\CU}{{\mathcal{U}}} \nc{\CV}{{\mathcal{V}}} \nc{\CW}{{\mathcal{W}}}
\nc{\CX}{{\mathcal{X}}} \nc{\CZ}{{\mathcal{Z}}}

\nc{\LT}{{\check{T}}} \nc{\hlambda}{{\hat{\lambda}}} \nc{\hulambda}{{\widehat{-\lambda-2\rho}}}
\nc{\hLambda}{{\widehat{\lambda+\Lambda}}} \nc{\hzeta}{{\hat{\zeta}}} \nc{\gen}{{\operatorname{gen}}}
\nc{\cM}{{\check{\mathcal M}}{}} \nc{\csM}{{\check{\mathcal A}}{}} \nc{\oM}{{\overset{\circ}{\mathcal M}}{}}
\nc{\obM}{{\overset{\circ}{\mathbf M}}{}} \nc{\oCA}{{\overset{\circ}{\mathcal A}}{}}
\nc{\obA}{{\overset{\circ}{\mathbf A}}{}} \nc{\ooM}{{\overset{\circ}{M}}{}}
\nc{\osM}{{\overset{\circ}{\mathsf M}}{}} \nc{\vM}{{\overset{\bullet}{\mathcal M}}{}}
\nc{\nM}{{\underset{\bullet}{\mathcal M}}{}} \nc{\oD}{{\overset{\circ}{\mathcal D}}{}}
\nc{\obD}{{\overset{\circ}{\mathbf D}}{}} \nc{\oA}{{\overset{\circ}{\mathbb A}}{}}
\nc{\op}{{\overset{\bullet}{\mathbf p}}{}} \nc{\cp}{{\overset{\circ}{\mathbf p}}{}}
\nc{\oU}{{\overset{\bullet}{\mathcal U}}{}} \nc{\oZ}{{\overset{\circ}{\mathcal Z}}{}}
\nc{\ofZ}{{\overset{\circ}{\mathfrak Z}}{}} \nc{\ovin}{{\overset{\circ}{\on{Vin}}}}

\nc{\fa}{{\mathfrak{a}}} \nc{\fb}{{\mathfrak{b}}} \nc{\fg}{{\mathfrak{g}}} \nc{\g}{{\mathfrak{g}}}
\nc{\fgl}{{\mathfrak{gl}}} \nc{\fh}{{\mathfrak{h}}} \nc{\fj}{{\mathfrak{j}}} \nc{\fm}{{\mathfrak{m}}}
\nc{\fn}{{\mathfrak{n}}} \nc{\fu}{{\mathfrak{u}}} \nc{\fp}{{\mathfrak{p}}} \nc{\frr}{{\mathfrak{r}}}
\nc{\fs}{{\mathfrak{s}}} \nc{\ft}{{\mathfrak{t}}} \nc{\fT}{{\mathfrak{T}}} \nc{\ofT}{{\overline{\mathfrak
T}}} \nc{\ofS}{{\overline{\mathfrak S}}} \nc{\fsl}{{\mathfrak{sl}}} \nc{\hsl}{{\widehat{\mathfrak{sl}}}}
\nc{\hgl}{{\widehat{\mathfrak{gl}}}} \nc{\hg}{{\widehat{\mathfrak{g}}}}
\nc{\chg}{{\widehat{\mathfrak{g}}}{}^\vee} \nc{\hn}{{\widehat{\mathfrak{n}}}}
\nc{\chn}{{\widehat{\mathfrak{n}}}{}^\vee}

\nc{\fA}{{\mathfrak{A}}} \nc{\fB}{{\mathfrak{B}}} \nc{\fC}{{\mathfrak{C}}} \nc{\fD}{{\mathfrak{D}}}
\nc{\fE}{{\mathfrak{E}}} \nc{\fF}{{\mathfrak{F}}} \nc{\fG}{{\mathfrak{G}}} \nc{\fI}{{\mathfrak{I}}}
\nc{\fK}{{\mathfrak{K}}} \nc{\fL}{{\mathfrak{L}}} \nc{\fM}{{\mathfrak{M}}} \nc{\fN}{{\mathfrak{N}}}
\nc{\frP}{{\mathfrak{P}}} \nc{\fS}{{\mathfrak S}} \nc{\fU}{{\mathfrak{U}}} \nc{\fX}{{\mathfrak{X}}}
\nc{\fY}{{\mathfrak{Y}}} \nc{\fZ}{{\mathfrak{Z}}}

\nc{\bb}{{\mathbf{b}}} \nc{\bc}{{\mathbf{c}}} \nc{\be}{{\mathbf{e}}} \nc{\bj}{{\mathbf{j}}}
\nc{\bn}{{\mathbf{n}}} \nc{\bp}{{\mathbf{p}}} \nc{\bq}{{\mathbf{q}}} \nc{\bfu}{{\mathbf{u}}}
\nc{\bv}{{\mathbf{v}}} \nc{\bx}{{\mathbf{x}}} \nc{\by}{{\mathbf{y}}} \nc{\bw}{{\mathbf{w}}}
\nc{\bA}{{\mathbf{A}}} \nc{\bB}{{\mathbf{B}}} \nc{\bC}{{\mathbf{C}}} \nc{\bF}{{\mathbf{F}}}
\nc{\bI}{{\mathbf{I}}} \nc{\bbI}{{\mathbf{I}}} \nc{\bK}{{\mathbf{K}}} \nc{\bD}{{\mathbf{D}}}
\nc{\bH}{{\mathbf{H}}} \nc{\bM}{{\mathbf{M}}} \nc{\bN}{{\mathbf{N}}} \nc{\bT}{{\mathbf{T}}}
\nc{\bV}{{\mathbf{V}}} \nc{\bW}{{\mathbf{W}}} \nc{\bX}{{\mathbf{X}}} \nc{\bP}{{\mathbf{P}}}
\nc{\bZ}{{\mathbf{Z}}} \nc{\bnu}{{\boldsymbol{\nu}}}

\nc{\sA}{{\mathsf{A}}} \nc{\sB}{{\mathsf{B}}} \nc{\sC}{{\mathsf{C}}} \nc{\sD}{{\mathsf{D}}}
\nc{\sF}{{\mathsf{F}}} \nc{\sK}{{\mathsf{K}}} \nc{\sM}{{\mathsf{M}}} \nc{\sO}{{\mathsf{O}}}
\nc{\sQ}{{\mathsf{Q}}} \nc{\sP}{{\mathsf{P}}} \nc{\sZ}{{\mathsf{Z}}} \nc{\sfp}{{\mathsf{p}}}
\nc{\sr}{{\mathsf{r}}} \nc{\sfb}{{\mathsf{b}}} \nc{\sfc}{{\mathsf{c}}} \nc{\sd}{{\mathsf{d}}}
\nc{\sfl}{{\mathsf{l}}}

\nc{\GG}{{\mathcal{G}}}

\nc{\BK}{{\bar{K}}}

\nc{\tB}{{\widetilde{\mathcal{B}}}} \nc{\tg}{{\widetilde{\mathfrak{g}}}}
\nc{\tG}{{\widetilde{G}}} \nc{\TM}{{\widetilde{\mathbb{M}}}{}} \nc{\tO}{{\widetilde{\mathbb{O}}}{}}
\nc{\tU}{{\widetilde{\mathfrak{U}}}{}} \nc{\TZ}{{\tilde{Z}}} \nc{\tx}{{\tilde{x}}} \nc{\tbv}{{\tilde{\bv}}}
\nc{\tfP}{{\widetilde{\mathfrak{P}}}{}} \nc{\tz}{{\tilde{\zeta}}} \nc{\tmu}{{\tilde{\mu}}}

\nc{\ul}{\underline} \nc{\ol}{\overline} \nc{\urho}{\underline{\rho}} \nc{\uB}{\underline{B}}
\nc{\uC}{{\underline{\mathbb{C}}}} \nc{\uc}{{\underline{c}}} \nc{\ucs}{{\underline{c},\operatorname{ss}}}
\nc{\ui}{\underline{i}} \nc{\uj}{\underline{j}} \nc{\ofP}{{\overline{\mathfrak{P}}}}
\nc{\oB}{{\overline{\mathcal{B}}}} \nc{\og}{{\overline{\mathfrak{g}}}} \nc{\oI}{{\overline{I}}}

\nc{\eps}{\varepsilon} \nc{\hrho}{{\hat{\rho}}}

\nc{\one}{{\mathbf{1}}} \nc{\two}{{\mathbf{t}}}

\nc{\Rep}{{\mathop{\operatorname{\rm Rep}}}} 
\nc{\Tot}{{\mathop{\operatorname{\rm Tot}}}} 
\nc{\Ker}{{\mathop{\operatorname{\rm Ker}}}} \nc{\Hilb}{{\mathop{\operatorname{\rm Hilb}}}}
\nc{\End}{{\mathop{\operatorname{\rm End}}}}
\nc{\Ext}{{\mathop{\operatorname{\rm Ext}}}} \nc{\Eext}{{{\mathcal{E}}xt}}
\nc{\Hom}{{\mathop{\operatorname{\rm Hom}}}} \nc{\RHom}{{\mathop{\operatorname{\rm RHom}}}}
\nc{\CHom}{{\mathop{\operatorname{{\mathcal{H}}\it om}}}} \nc{\GL}{{\mathop{\operatorname{\rm GL}}}}
\nc{\gr}{{\mathop{\operatorname{\rm gr}}}} \nc{\Id}{{\mathop{\operatorname{\rm Id}}}}
\nc{\Ind}{{\mathop{\operatorname{\rm Ind}}}} \nc{\defi}{{\mathop{\operatorname{\rm def}}}}
\nc{\length}{{\mathop{\operatorname{\rm length}}}} \nc{\supp}{{\mathop{\operatorname{\rm supp}}}}

\nc{\Cliff}{{\mathsf{Cliff}}} 
\nc{\Fl}{{\mathsf{Fl}}}
\nc{\Fib}{{\mathsf{Fib}}} \nc{\Coh}{{\mathsf{Coh}}} \nc{\FCoh}{{\mathsf{FCoh}}}

\nc{\reg}{{\text{\rm reg}}}

\nc{\cplus}{{\mathbf{C}_+}} \nc{\cminus}{{\mathbf{C}_-}} \nc{\cthree}{{\mathbf{C}_*}} \nc{\Qbar}{{\bar{Q}}}

\nc{\bh}{{\bar{h}}} \nc{\bOmega}{{\overline{\Omega}}}

\nc{\seq}[1]{\stackrel{#1}{\sim}}

\newcommand{\beq}{\begin{equation}} \newcommand{\eeq}{\end{equation}} \renewcommand{\proof}{{\it Proof }}
\newcommand{\proofpt}{{\it Proof. }}

\newcommand{\imbed}{\hookrightarrow}

\newcommand{\oplusl}{\bigoplus\limits} \newcommand{\cupl}{\bigcup\limits}

\newcommand{\la}{\lambda}
\newcommand{\La}{\Lambda}

\newcommand{\Z}{{\mathcal{Z}}} 

\nc{\hatt}{\widehat}

\author{Roman Bezrukavnikov, Michael Finkelberg and Victor Ostrik} \title[Character $D$-modules via Drinfeld
center]{Character $D$-modules via Drinfeld center of Harish-Chandra bimodules}

\begin{document}

\begin{abstract} The category of character D-modules is realized as Drinfeld center of the abelian monoidal
category of Harish-Chandra bimodules. Tensor product of Harish-Chandra bimodules
 is related to convolution
of D-modules via the long intertwining functor (Radon transform) by a result of \cite{BG}. Exactness
property of the long intertwining functor on a cell subquotient of the Harish-Chandra bimodules category
shows that the
 truncated convolution category of \cite{Ltc} can be realized as a subquotient
 of the category of Harish-Chandra bimodules. Together with the description
 of the truncated convolution category \cite{BFO} this
allows us to derive (under a mild technical assumption)
a classification of irreducible character sheaves over $\Ce$ obtained by Lusztig by a different method.

We also
  give a simple description for the top cohomology of
 convolution of character sheaves  over $\Ce$ in a given cell
modulo smaller cells and relate the so-called Harish-Chandra functor to Verdier specialization in the De
Concini-Procesi compactification.
 \end{abstract}

\maketitle 
\tableofcontents

\section{Introduction} The goal of this paper is to reach a better understanding of some aspects of
Lusztig's theory of character sheaves.

Let $G$ be a semi-simple split algebraic group with a Borel subgroup $B$. Character sheaves were developed
in~\cite{L85} as a tool for description of characters of the finite group $G(\Fq)$. For simplicity, in the
introduction we restrict attention to the case of unipotent character sheaves.

The set $Rep_{unip}(G(\Fq))$ of unipotent irreducible representations of $G(\Fq)$ contains the set
$PS(G(\Fq))$ of representations generated by a Borel invariant vector (also known as the principal series).
It is well known that $PS(G(\Fq))$ is in bijection with simple modules over the Hecke algebra
$H=\Ce[B(\Fq)\backslash G(\Fq)/B(\Fq)]$ which is a semi-simple algebra isomorphic to $\Ce[W]$. Thus
the cardinality of $PS(G(\Fq))$ equals the dimension of the center $Z(H)$.

Though for  $G$  of type $A_n$ we have $Rep_{unip}(G(\Fq))=PS(G(\Fq))$, in general $PS(G(\Fq))\subsetneq
Rep_{unip}(G(\Fq))$. The main point of the present paper is related to the fact that the cardinality of
$Rep_{unip}(G(\Fq))$ equals the rank of the Grothendieck group of the {\em categorical (Drinfeld) center of
the categorical Hecke algebra}. In fact,
this categorical center is identified with the category of character sheaves. Furthermore, the irreducible
character sheaves correspond to irreducible objects in the categorical center of the {\em categorical
asymptotic Hecke algebra}.

Let us now make this more precise. By the categorical Hecke algebra we mean, following a well established
pattern, the  $B$ equivariant  derived category of sheaves on the flag space $G/B$ (or a technical
variation, such as the derived category of $T\times T$ unipotently monodromic $N$ equivariant sheaves on
$G/N$; we will ignore the difference in the rest of the Introduction). In fact, we work on the De Rham side
of the Riemann-Hilbert correspondence over a field of characteristic zero, thus the original motivating
setting of groups over a finite field will not be mentioned below (we plan to return to it in a future
publication).

Instead we consider the category of $B$-equivariant $D$-modules on $G/B$. This can also be realized as the
category of $G$-equivariant $D$-modules on $G/B\times G/B$. By the Beilinson-Bernstein Localization Theorem this
category is equivalent to the category $\Har^0$ of Harish-Chandra bimodules with trivial central character
(more general central characters are considered in the main body of the paper; also, for technical reasons
it is better to work with the category $\Har^{\hat{0}}$ of Harish-Chandra bimodules with a trivial generalized central character).  The tensor product of bimodules provides $\Har^{\hat{0}}$ with a monoidal structure, taking the left derived functor of the right exact tensor product we get a monoidal triangulated category $(D\Har^{\hat{0}}, \Lotimes _U)$.

The natural exact equivalence of $\Har^{\hat{0}}$ with a category of $\CalD$-modules does {\em not}
intertwine this monoidal structure with convolution of complexes of $\CalD$-modules. However, the
composition of this equivalence with the long intertwining functor $\bI_{w_0}$ turns out to be monoidal, by
a result of \cite{BG}. Thus $(D\Har^{\hat{0}}, \Lotimes _U)$ can also be thought of as a categorification\footnote{See e.g.
\cite{Str1} for a somewhat different (though related) approach to categorifying the Hecke algebra and its modules.} of
the Hecke algebra $H$.

Our first result, explained in section 3, asserts that the Drinfeld center of the abelian category
$\Har^{\hat{0}}$ is canonically equivalent
 to the category of character $\CalD$-modules on the group $G$. An
essential step in the argument is an exactness property of the functor $\HC$, which is a categorification of the
action of the group algebra of $G(\Fq)$ on a principal series module. The functor $\HC$ itself is not
exact, however, we show that its composition with an appropriate direct summand in the functor of global
sections is exact. Using a result of Beilinson-Bernstein \cite{BB}, the latter statement can be reformulated
geometrically, as exactness of the composition of $\HC$ with the (inverse to the)  long intertwining functor.\footnote{What we call an intertwining functor is sometimes called, depending on the context and on the author, a Radon transform functor, or a shuffling functor, or a twisting functor.} In
the last section 6 we sketch a proof of the fact that this exact functor can be described as Verdier
specialization functor from $D$-modules on $G$ to $D$-modules on the punctured normal bundle to a stratum in the De Concini -- Procesi compactification (Corollary~\ref{spec}) -- this statement allows one to reprove the exactness property, but is also of independent interest.

In section  4 we 
consider a  filtration 
on the category of Harish-Chandra bimodules according to the support of a bimodule.
It turns out that the subcategory of semi-simple objects in the subquotient category is closed under tensor product;
the resulting semi-simple monoidal category turns
out to coincide with the truncated convolution category introduced in \cite{Ltc}. This follows from the fact that the functor on the subquotient category induced by $\bI_{w_0}$ is exact up to a shift.
The relation between Lusztig's asymptotic Hecke algebra
 $J$ (which is well-known to be
identified with the Grothendieck group of the truncated convolution category) and tensor product of
Harish-Chandra bimodules is implicit in the work of A.~Joseph, cf. e.g. \cite[2.3(i)]{J1}.

In section  5 we check that irreducible objects in the center of the category of Harish-Chandra bimodules
are in bijection with irreducible objects in the center of the category of semi-simple objects in the
associated graded category with respect to the support filtration -- i.e., roughly speaking, taking center of the monoidal category commutes with passing to the
 semisimplification of the associated graded. The argument uses a technical statement that the action of the center of Harish-Chandra modules on the subquotient category preserves the subcategory of semi-simple objects. This is deduced from a recent result of I.~Losev \cite{Los} recalled in \ref{recL}.

This is used to deduce Lusztig's classification of irreducible character sheaves.
Namely, to each 2-sided cell $\uc$ in $W$ Lusztig assigned a finite group ${\mathcal G}_\uc$.
He has shown that the set of irreducible unipotent character sheaves is in bijection with the set $\cupl_{\uc} \{ (\gamma, \psi)\}/{\mathcal G}_\uc$ where $\gamma\in {\mathcal G}_\uc$ and $\psi$ is
an irreducible representation of the centralizer of $\gamma$ in ${\mathcal G}_\uc$, i.e. by the union over $\uc$ of the sets of irreducible conjugation equivariant sheaves
on ${\mathcal G}_{\uc}$.
The description of the truncated convolution category from \cite{BFO} (conjectured
by Lusztig and recalled below before Corollary \ref{cor54}) together with 2-Morita invariance of categorical center implies that irreducible objects in the center of truncated convolution category are indexed by this set.
Since the semi-simple part of the associated graded
 for the support filtration on Harish-Chandra bimodules is identified with the
truncated convolution category, we get the classification result.

This also implies that truncated convolution of character sheaves modulo smaller cells
corresponds to convolution of conjugation equivariant sheaves on
the finite group ${\mathcal G}_{\uc}$ (though the definition of the braided category of sheaves on ${\mathcal G}_{\uc}$ has to be modified in some cases
 in order to upgrade this to a tensor equivalence, see Corollary \ref{cor54}).
We note in passing that a different statement of a similar nature  was conjectured by
Lusztig~\cite{L04}. While in the present paper we define truncated convolution of character sheaves as the
component of their convolution  concentrated in a certain homological
 degree, Lusztig in {\em loc. cit.} considers a component of a certain weight.

Let us remark that the description of the truncated convolution categories
which plays a crucial role in our proof of the classification was proved in \cite{BFO}
in a somewhat indirect way, using, in particular, the embedding of finite Hecke
algebra into the affine one, and its categorification. A more direct approach based on representations of finite $W$-algebras will be presented in a forthcoming paper \cite{LO}, see footnote after Theorem \ref{december} below.

\medskip

Let us mention also that the idea to relate Drinfeld center to character sheaves has been implemented in a different way by Ben-Zvi and Nadler \cite{BZN}, who have obtained a result resembling our Theorem \ref{davids}. However, their approach is different from ours in that they work with the full triangulated category of equivariant complexes
 on $G/B$ rather than with the abelian heart of a particular $t$-structure.
The categorical center construction does not produce reasonable results when applied to triangulated categories, and the recent formalism of derived algebraic geometry
is employed in \cite{BZN} to overcome this difficulty. Our observation on the exactness property of the horocycle transform (Proposition \ref{HotKas}(b)) allows us to bypass that difficulty by working with abelian categories; thus the technique of the present paper is more elementary than that of \cite{BZN}. Also, the application to classification of irreducible character sheaves, which was one of our main motivations, does not appear in \cite{BZN} and is unlikely to be accessible by its methods alone.

\medskip

The statements established in this paper for $D$-modules have obvious analogues in the setting of perverse sheaves
in classical topology over $\Ce$ and $l$-adic sheaves. Although the proofs in the first setting could be
deduced from the results below via
 Riemann-Hilbert equivalence, it would be desirable to find direct
geometric proofs applicable also to $l$-adic sheaves. We plan to return to this problem in a future publication.

{\bf Acknowledgements.} We thank A.~Joseph and D.~Vogan for help with references. We are grateful to
V.~Ginzburg, I.~Losev and C.~Stroppel
 for useful discussions, and to I.~Mirkovi\' c for sharing some related ideas (see Remark \ref{mir}).
 The first author thanks A.~Beilinson,
V.~Drinfeld and D.~Vogan for the stimulating opportunity to present the preliminary version of the results
in a seminar. R.B. was partially supported by the DARPA grant  HR0011-04-1-0031 and NSF grants DMS-0625234,
DMS-0854764;
M.F. was partially supported by the AG Laboratory HSE, RF government
grant, ag. 11.G34.31.0023, RFBR grant 09-01-00242, the Ministry of
Education and Science of Russian Federation grant No.
2010-1.3.1-111-017-029, and the Science Foundation of the NRU-HSE
award 11-09-0033; and V.O. was partially supported by the
NSF grant DMS-0602263.

\section{Setup}






\subsection{Notations related to $G$} Let $G$ be an almost simple complex algebraic group. We denote by $T$
the abstract Cartan group of $G$. We also fix a Cartan and Borel subgroup $T\subset B\subset G$. The Weyl
group of $G$ is denoted by $W$. We denote by $N$ the unipotent radical of $B$. The Lie algebras of $T\subset G$ are denoted by $\ft\subset\fg$. We denote by $\Lambda$ the weight lattice of $T$. We denote by
$2\rho\in\ft^*$ the sum of positive roots of $T\subset B$. We fix a dominant regular 
rational weight
$\lambda\in\ft^*$, that is $\lambda$ does not lie on any coroot hyperplane shifted by $-\rho$, and the value of $\lambda$ at any positive coroot is not a negative integer. We consider $\zeta:=\exp(2\pi i\lambda)$ as an element of the dual torus $\LT$. We denote by $\CB$ the flag variety $G/B$, and we denote by $\tB$ the base affine space $G/N$. The horocycle space $\CY$ is the quotient of $G/N\times G/N$ by the diagonal right action of Cartan $T$. It is a $G$-equivariant $T$-torsor over $\CB\times\CB$.

\subsection{Sheaves on horocycle space} We consider the category $\CP^\zeta$ (resp. $\CP^\hzeta$) of
$G$-equivariant $T$-monodromic perverse sheaves on $\CY$ with monodromy (resp. generalized monodromy)
$\zeta$. We denote by $D\CP^\zeta$ (resp. $D\CP^\hzeta$) the $G$-equivariant constructible
$\zeta$-monodromic (resp. generalized $\zeta$-monodromic) derived category on $\CY$. By the Riemann-Hilbert
correspondence, $\CP^\hzeta$ is equivalent to a certain category $\CM^\hLambda$ of $G$-equivariant
$\CalD_\CY$-modules. Let $\on{p}$ stand for the projection $\tB\times\tB\to\CY$. Then the functor
$\on{p}^![-\dim T]$ identifies $\CM^\hLambda$ with the category of $G$-equivariant, $T$-equivariant (with
respect to the diagonal right action) $\CalD_{\tB\times\tB}$-modules such that $\BC[\ft^*]\otimes\BC[\ft^*]$
(acting as infinitesimal translations along the right action of $T\times T$) acts locally finitely with
generalized eigenvalues in $(\lambda+\Lambda,-\lambda+\Lambda)$. We denote by $D\CM^\hLambda$ the
$G$-equivariant derived category of $\CM^\hLambda$.

For the rational weights $\mu,\nu$ we denote by $\CM^{\widehat{\mu+\Lambda},\widehat{\nu+\Lambda}}$ the
category of $G$-equivariant $\CalD_{\tB\times\tB}$ modules ($T$-equivariance is no longer required) such that
$\BC[\ft^*]\otimes\BC[\ft^*]$ acts locally finitely with generalized eigenvalues in
$(\mu+\Lambda,\nu+\Lambda)$. We denote by $D\CM^{\widehat{\mu+\Lambda},\widehat{\nu+\Lambda}}$ the
$G$-equivariant derived category of $\CM^{\widehat{\mu+\Lambda},\widehat{\nu+\Lambda}}$. The functor
$\on{p}^o:=\on{p}^![-\dim T]:\ D\CM^\hLambda\to D\CM^{\widehat{\lambda+\Lambda},\widehat{-\lambda+\Lambda}}$
identifies $\CM^\hLambda$ with a full subcategory of
$\CM^{\widehat{\lambda+\Lambda},\widehat{-\lambda+\Lambda}}$.

\subsection{Harish-Chandra bimodules} We denote by $U=U(\fg)$ the universal enveloping algebra of $\fg$. We
denote by $U_\lambda$ the quotient of $U$ by the ideal generated by the maximal ideal $Z_\lambda$ of the
Harish-Chandra center $Z(U)$, that is the ideal corresponding to the maximal ideal $I_\lambda\in\ft$ under the Harish-Chandra isomorphism. We denote by $U_\hlambda$ the completion of $U$ at this ideal. We have
$U_\hlambda^{opp}\simeq U_\hulambda=U_{\widehat{-w_0\lambda}}$ where $w_0\in W$ is the longest element.
Furthermore, we consider an inductive system of ideals of $Z(U)$ of finite codimension supported at
$\lambda+\Lambda$. It gives rise to the projective system of quotients of $U$ by the two-sided ideals
generated by the above ideals of $Z(U)$. The projective limit of this system is denoted by $U_\hLambda$.
Note that we have an embedding $U_\hlambda\hookrightarrow U_\hLambda$.

We denote by $\Har^\lambda$ the category of finitely generated (as a left module; equivalently, as a right
module) Harish-Chandra bimodules over $U_\lambda$, that is, $U_\lambda$-bimodules such that the adjoint
$\fg$-action is locally finite. Note that a $U_\lambda$-bimodule is the same as a $U_\lambda\otimes
U_\lambda^{opp}= U_\lambda\otimes U_{-\lambda-2\rho}$-module, and hence a Harish-Chandra
$U_\lambda$-bimodule is the same as a $U_\lambda\otimes U_{-\lambda-2\rho}$-module such that the diagonal
action of $\fg$ is locally finite. We denote by $\Har^\hlambda$ the category of (finitely generated as a
left $U$-module) Harish-Chandra bimodules over $U_\hlambda$. Note that a $U_\hlambda$-bimodule is the same
as a $U_\hlambda\otimes U_\hlambda^{opp}= U_\hlambda\otimes U_\hulambda$-module. We denote by
$D\Har^\lambda$ (resp. $D\Har^\hlambda$) the bounded derived category of $\Har^\lambda$ (resp.
$\Har^\hlambda$).

For $M\in\CM^{\widehat{\mu'+\Lambda},\widehat{\nu'+\Lambda}}$, the action of $\BC[\ft^*]\otimes\BC[\ft^*]$
on the global sections $\Gamma(\tB\times\tB,M)$ is locally finite, and for $\mu\in\mu'+\Lambda,\
\nu\in\nu'+\Lambda$, we denote by $\Gamma^{\hat{\mu},\hat{\nu}}(M)$ the maximal summand where
$\BC[\ft^*]\otimes\BC[\ft^*]$ acts with the generalized eigenvalue $(\mu,\nu)$. It extends to the derived
category $D\CM^{\widehat{\mu'+\Lambda},\widehat{\nu'+\Lambda}}$, and for $\mu=\lambda,\ \nu=-\lambda-2\rho$,
gives rise to an equivalence of categories $R\Gamma^{\hlambda,\hulambda}:\
D\CM^{\hLambda,\widehat{-\lambda+\Lambda}} \iso D\Har^\hlambda$.

\subsection{Convolution} \label{convo} We have the convolution functor $*:\
D\CM^{\hLambda,\widehat{-\lambda+\Lambda}}\times D\CM^{\hLambda,\widehat{-\lambda+\Lambda}}\to
D\CM^{\hLambda,\widehat{-\lambda+\Lambda}}$, cf.~\cite{BG}~(5.6),~(5.13). If $p$ stands for the projection
$\tB\times\tB\times\tB\to\tB\times\tB$ along the middle factor, and $\Delta$ stands for the inclusion
$\tB\times\tB\times\tB\hookrightarrow\tB\times\tB\times\tB\times\tB,\ (x,y,z)\mapsto(x,y,y,z)$, then
$M_1*M_2=p_*\Delta^![\dim\tB] (M_1\boxtimes M_2)$. Let us stress that $p_*\ne p_!$ since $\tB$ is not
proper.

We also have the convolution functor $*:\ D\Har^\lambda\times D\Har^\lambda\to D\Har^\lambda$ (resp.
$D\Har^\hlambda\times D\Har^\hlambda\to D\Har^\hlambda$), see {\em loc. cit}. Note that the convolution of
Harish-Chandra bimodules is the left derived functor of the right exact bifunctor $(A,B)\mapsto
A\otimes_{U_\hlambda}B,\ \Har^\hlambda\times\Har^\hlambda\to\Har^\hlambda$. Finally, for arbitrary
$\CalD$-modules $C$ on $G$, and $M$ on $\tB\times\tB$ we set $C*M:=a_*(C\boxtimes M)$ where $a:\
G\times\tB\times\tB\to\tB\times\tB$ is the action morphism.

\begin{lem}\label{conv_abs} Let $X$ be an algebraic variety with an action of an (affine) algebraic group
$G$. Then for $A\in D^b(\CalD$-mod$(G))$, $M\in D^b(\CalD$-mod$(X))$ we have a canonical isomorphism
$$R\Gamma(A*M)\cong R\Gamma(A)\Lotimes_U R\Gamma(M).$$ \end{lem} \proofpt Replacing $A$, $M$ by appropriate
resolutions we reduce to the case when $A$ is a locally free $\CalD$-module and $M=j_*(L)$ where $j:\ Y\to
X$ is an embedding of an open affine subvariety. Let $a:\ G\times Y\to X$ be the action map and let $Vert$
denote the space of vertical vector fields for that map.
 Then the left hand side lies in homological degree zero and equals the module
of coinvariants of $Vert$ acting on $\Gamma((A\boxtimes L)\otimes\omega_{G\times U} \otimes \omega_Y^{-1})$,
where $\omega$ denotes the line bundle of top degree forms. The vector fields coming from the $G$-action on
$G\times X$, $g:(g_1,x)\mapsto (g_1g^{-1},gx)$ generate the vector bundle $Vert$. Moreover, the bi-invariant
volume form on $G$ yields a canonical $G$-invariant section of the line bundle $\omega_{G\times Y} \otimes
\omega_Y^{-1}$. This yields the desired isomorphism. \qed

\begin{cor}\label{11} \label{hse} a) For $\CalD$-modules $A$ on $G$, and $M$ on $\tB\times\tB$ we have
$R\Gamma(\tB\times\tB,A*M)=R\Gamma(M)\otimes_U\Gamma(G,A)$.

b) For $M_1,M_2\in D\CM^{\hLambda,\widehat{-\lambda+\Lambda}}$ we have a natural isomorphism
\begin{equation} \label{bg} R\Gamma^{\hlambda,\hulambda}(M_1*M_2)\simeq
R\Gamma^{\hlambda,\hulambda}(M_1)*R\Gamma^{\hlambda,\hulambda}(M_2), \end{equation} \end{cor}

\proofpt a) is immediate from the previous Lemma.

(b)  is proved in~\cite{BG},~Proposition~5.11. \qed

\subsection{Character sheaves}\label{CSdef} We consider the following diagram of $G$-varieties and $G$-equivariant maps:
$$G\stackrel{pr}{\longleftarrow}G\times\CB\stackrel{f}{\longrightarrow}\CY.$$ In this diagram, the map $pr$
is given by $pr(g,x):=g$. To define the map $f$, we think of $\CB$ as $\tB/T$, and for a representative
${\tilde x}\in\tB$ of $x\in\CB$ we set $f(g,x):=(g{\tilde x},{\tilde x})$. Following~\cite{G}, we consider
the functor $\CH=pr_*f^*[\dim\CB]:\ D\CP^\hzeta\to D_G(G)$ to the $G$-equivariant derived constructible
category on $G$. The minimal triangulated full subcategory of $D_G(G)$ containing all the objects of the
form $\CH(A),\ A\in D\CP^\hzeta$, will be denoted by $D\Char^\hzeta(G)$. The intersection of
$D\Char^\hzeta(G)$ with the abelian subcategory of perverse sheaves on $G$ will be denoted by
$\Char^\hzeta$. The full subcategory of semisimple objects of $\Char^\hzeta$ will be denoted $\Char^\zeta$
--- these are direct sums of Lusztig's character sheaves with monodromy $\zeta$. The functor $\CH$ admits
the right adjoint functor $\HC:=f_*pr^![-\dim\CB]:\ D\Char^\hzeta(G)\to D\CP^\hzeta$.

Let $\CalD(G)$ stand for the ring of differential operators on $G$. We have $\CalD(G)=U\ltimes\CO(G)$. 

 We denote by $\Char^\hLambda$ the category of (finitely generated over $\CalD(G)$)
$G$-equivariant $\CalD(G)$-modules such that the action of
$Z(U)\subset U\subset\CalD(G)$ embedded as left invariant differential operators is locally finite and has
generalized eigenvalues in $\lambda+\Lambda$ (equivalently, such that the action of $Z(U)\subset
U\subset\CalD(G)$ embedded as right invariant differential operators has generalized eigenvalues in
$-\lambda-2\rho+\Lambda$). Such $\CalD(G)$-modules are automatically holonomic regular (see~\cite{G}). We
denote by $D\Char^\hLambda$ the $G$-equivariant derived category of $\CalD_\hLambda(G)$-modules. The
Riemann-Hilbert correspondence gives rise to an equivalence $\Char^\hzeta\simeq\Char^\hLambda$ (resp.
$D\Char^\hzeta\simeq D\Char^\hLambda$), see {\em loc. cit}. Clearly the action
of $\CalD(G)$ on  $M\in\Char^\hLambda$
factors through  $U_\hLambda\ltimes\CO(G)$, where
$U_\hLambda\ltimes\CO(G)$ denotes
the subalgebra of $ad(\g)$-finite vectors 
 in the inverse limit of images of  $\CalD(G)$ in endomorphisms of modules in $\Char^\hLambda$,
 this is a topological ring isomorphic to $U_\hLambda\ltimes\CO(G)$ as a vector space.
 We can identify $\Char^\hLambda$ with the category of $G$-equivariant discrete
 finitely generated $U_\hLambda\ltimes\CO(G)$-modules.

We will keep the names $\CH:\ D\CM^\hLambda\to D\Char^\hLambda;\ \HC:\ D\Char^\hLambda\to D\CM^\hLambda$ for
the functors going to the same named functors under the Riemann-Hilbert correspondence.

\section{Character sheaves as Drinfeld center of Harish-Chandra bimodules}

\subsection{Exactness of $\HC$}
\label{huji}
For $A\in D\Char^\hLambda$ we define $\Gamma^\hlambda(G,A)$ as follows. Note that $U_\hLambda$ is the direct
product $\prod_{\mu\in(\lambda+\Lambda)/W_\zeta}U_{\hat\mu}$ over the set of representatives of
$W_\zeta$-orbits (with respect to the dot action) in $\lambda+\Lambda$. Here $W_\zeta\subset W$ is the
stabilizer of $\zeta\in\LT$ in $W$. Equivalently, $W_\zeta$ is the stabilizer of $\lambda+\Lambda$ in
$\ft^*/\Lambda$. Hence the restriction of a $\CalD_\hLambda(G)$-module $A$ to $U_\hLambda$ has a maximal
direct summand $\Gamma^\hlambda(G,A)$ such that the action of any $U_{\hat\mu},\ \mu\not\in
W_\zeta\cdot\lambda$, on $\Gamma^\hlambda(G,A)$ is zero. Here we think of $U$ as acting by the left
invariant differential operators. Equivalently, $\Gamma^\hlambda(G,A)$ is the maximal direct summand of $A$
such that the action of any $U_{\hat\mu},\ \mu\not\in W_\zeta\cdot(-\lambda-2\rho)$, as right invariant
differential operators, is zero.

\begin{prop} \label{HotKas} a) For $A\in D\Char^\hLambda$ we have \begin{equation} \label{KaHo}
R\Gamma^{\hlambda,\hulambda}(\on{p}^o\HC(A))=\Gamma^\hlambda(G,A) \end{equation}

b) The functor $R\Gamma^{\hlambda,\hulambda}\circ\on{p}^o\circ\HC:\ D\Char^\hLambda\to D\Har^\hlambda$ is
exact, i.e. takes $\Char^\hLambda$ to $\Har^\hlambda$. \end{prop}

\proofpt Clearly, the functor $\Gamma^\hlambda(G,-)$ is exact, so it remains to prove~(\ref{KaHo}).

The preimage of diagonal  $\Delta_\CB\subset \CB\times\CB$ in $\CY$ is isomorphic to $\CB\times T$. Let
$\delta$ denote the embedding of $\CB\times \{1\}$ into $\CY$. It is easy to see that $\HC(A)\cong
A*\delta_*(\CO)$ (notations of subsection \ref{convo}). On the other hand, it is easy to see that
$R\Gamma(\delta_*(\CO)\otimes pr_2^*\omega)\cong \tilde U$ where $\tilde U=U\otimes _{S(\ft)^W} S(\ft)$
(here $S(\ft)^W$ is identified with the center of $U$ by means of Harish-Chandra isomorphism and $pr_2$
denotes the second projection $\CY\to \CB$). Thus \eqref{KaHo} follows from Corollary \ref{11}(a).
 \qed

For future use we also record an algebraic description of the functor $\CH$.
Recall that the ring $\CalD_\hLambda(G)$ contains
$U_\hlambda\subset U_\hLambda$ as the left invariant differential operators. It also contains $U_\hulambda$
as the right invariant differential operators, and the images of $U_\hlambda$ and $U_\hulambda$ commute.
Thus $\CalD_\hLambda(G)$ has a structure of a right module over $U_\hlambda\otimes U_\hulambda$.

\begin{cor}\label{corHoKa}
For a regular weight $\lambda$ and $M\in D\CM^\hLambda$ we have a canonical isomorphism:
\begin{equation}
\label{HoKa} R\Gamma(G,\CH(M))= \CalD_\hLambda(G)\stackrel{L}{\otimes}_{U_\hlambda\otimes U_\hulambda}
R\Gamma^{\hlambda,\hulambda}(\on{p}^oM) \end{equation}
\end{cor}

\proofpt The standard adjunction between induction and restriction implies that the functor
in the right hand side of \eqref{HoKa} is left adjoint to $\Gamma^\hlambda$. Since
$\CH$ is right adjoint to $\HC$, we get the desired isomorphism. \qed

\begin{rem}  The isomorphism \eqref{HoKa} is a straightforward generalization
of a result of Hotta--Kashiwara \cite[Theorem 1]{HK}. Thus alternatively we could have
deduced \eqref{HoKa} from \cite{HK} and \eqref{KaHo} from \eqref{HoKa}.
 \end{rem}

\subsection{The twisted Harish-Chandra functor}\label{inter} In this subsection
 we present a geometric reformulation of Proposition \ref{HotKas}(b).

Recall the definition of an {\em intertwining functor} $\bI_{w_0}$, cf.~\cite{BB}\footnote{Note that our
notations do not agree with that of~\cite{BB}, our $\bI_{w_0}$ is denoted by $\bI_{w_0}^{-1}$ in {\em loc.
cit.}} %
and~\cite{BG},\footnote{Our definition of an interwining functor is different from that of
\cite{BG}, though one can show that the two definitions
 are equivalent. By our definition
$\bI_{w_0}$ is a (shriek) convolution with a certain $G$-equivariant $D$-module on $\tB^2$. In \cite{BG} the same functor is described as a convolution with a certain pro-object in the category of $G$-equivariant
$D$-modules on $\tB^2$ which are also {\em monodromic} with respect to $T^2$.}
 pp.~18-19. Let $U\subset \tB^2$ be the free $G\times T$ orbit,
and let $V\subset \tB^4$ be given by $V=\{(x_1,x_2,y_1,y_2\ |\ x_1=y_1, \, (x_2,y_2)\in U\}$. Then the
intertwining functor is given by $\bI_{w_0}=pr_{y!}pr_x^*$, where $pr_x,\, pr_y:V\to \tB^2$ send
$(x_1,x_2,y_1,y_2)$ to $(x_1,x_2)$, $(y_1,y_2)$ respectively.

It is a standard fact that the $\bI_{w_0}$ restricts to an equivalence $\bI_{w_0}: \
D\CM^{\hLambda,\widehat{-\lambda+\Lambda}} \to D\CM^{\hLambda,\widehat{-w_0\lambda+\Lambda}}$. The inverse
equivalence is $\bI_{w_0}^{-1}=pr_{y*}pr_x^!$.


Set $\CX:=(G/N^+\times G/N^-)/T$. Here $N^+=N$ is the unipotent radical of $B$, while $N^-$ is the unipotent radical of the opposite Borel subgroup $B^-$. Equivalently, we consider the following embedding $T\hookrightarrow T\times T,\ t\mapsto(t,w_0t)$. Then $\CX$ is the quotient of $\tB\times\tB$ with respect to the action of $T$ embedded as above into $T\times T$ acting on the right. It is equipped with the action of $T$ (embedded into $T\times T$ by $t\mapsto(t,1)$). We consider the $G$-equivariant $T$-monodromic derived category $D\CM^{\widehat{\lambda+\Lambda}}_{w_0}$ of $\CalD_\CX$-modules (here $\lambda+\Lambda$ is the generalized eigenvalue of the right action of the first copy of $\ft$)). This category can be realized as a subcategory of $D\CM^{\hLambda,\widehat{-w_0\lambda+\Lambda}}$; namely $\CF\mapsto q^![-\dim T]\CF$ where $q:\ \tB\times\tB\to\CX$ is the natural projection. Recall that $D\CM^{\widehat{\lambda+\Lambda}}$ is
realized as a subcategory of $D\CM^{\widehat{\lambda+\Lambda},\widehat{-\lambda+\Lambda}}$.

The intertwining functor $\bI_{w_0}$ gives rise to the same named equivalence $D\CM^{\hLambda}\iso
D\CM^{\widehat{-w_0+\Lambda}}_{w_0}$.

\begin{cor}\label{exa} The functor $\bI_{w_0} \circ \HC $ is exact. \end{cor}

\proofpt Since both $\lambda$ and $-w_0\lambda$ are dominant regular, the functor $\Gamma^{\hat \lambda,
\widehat{-w_0(\lambda)}}$ is exact and faithful by Localization Theorem. Thus it suffices to show that
$R\Gamma^{\hat \lambda, \widehat{-w_0(\lambda)}} \circ \bI_{w_0} \circ \HC$ is exact. By (a straightforward
generalization of) the result of \cite{BB} this functor is isomorphic to
$R\Gamma^{\hlambda,\widehat{-\lambda-2\rho}} \circ \HC$, thus it is exact by Proposition \ref{HotKas}(b).
\qed

\begin{rem} In fact, the exact functor described in the Corollary can be expressed via
 {\em Verdier specialization}.  This stronger statement is proved in the last
section 6. \end{rem}

\subsection{The Drinfeld center of $\Har^\hlambda$}\label{Drcsec}
 For a monoidal category $\CA$ we denote by $\CZ(\CA)$
the Drinfeld center of $\CA$ (see e.g.~\cite{Ka},~XIII.4). The functor $\on{p}^o\circ\HC:\
D\Char^\hLambda\to D\CM^{\hLambda,\widehat{-\lambda+\Lambda}}$ is equipped with a canonical central
structure (with respect to the convolution $*$) which gives rise to the same named functor $\HC:\
D\Char^\hLambda\to \CZ(D\CM^{\hLambda,\widehat{-\lambda+\Lambda}})$. In effect, for a $\CalD$-module $A$ on
$G$, and $M$ on $\tB\times\tB$ let us denote by $A*_1M$ (resp. $A*_2M$) the action along the first (resp.
second) copy of $\tB$. We have the canonical isomorphisms $A*_1M\simeq\HC(A)*M,\ A*_2M\simeq M*\HC(A)$.
Also, for a conjugation-equivariant $A$, and $M$ equivariant with respect to the diagonal action of $G$, we
have $A*_1M\simeq A*_2M$. The composition of the above isomorphisms provides $\HC$ with a central structure.
Since $R\Gamma^{\hlambda,\hulambda}:\ D\CM^{\hLambda,\widehat{-\lambda+\Lambda}}\to D\Har^\hlambda$ is a
monoidal equivalence, we obtain the functor $R\Gamma^{\hlambda,\hulambda}\circ\on{p}^o\circ\HC:\
D\Char^\hLambda\to\CZ(D\Har^\hlambda)$. Recall that the category of Harish-Chandra bimodules $\Har^\hlambda$
is equipped with the monoidal structure $(A,B)\mapsto A\otimes_{U_\hlambda}B$ whose left derived bifunctor
is the convolution $*$ on $D\Har^\hlambda$. Now Proposition~\ref{HotKas} implies that the exact functor
$R\Gamma^{\hlambda,\hulambda}\circ\on{p}^o\circ\HC$ from the {\em abelian} category $\Char^\hLambda$ to the
{\em abelian} monoidal category $\Har^\hlambda$ has a canonical central structure which gives rise to the
exact functor $\hc:\ \Char^\hLambda\to\CZ(\Har^\hlambda)$.

The rest of the section is devoted to the proof of the following

\begin{thm} \label{davids} $\hc:\ \Char^\hLambda\to\CZ(\Har^\hlambda)$ is an equivalence of abelian
categories. \end{thm}

\proofpt Since any irreducible object $A$ of $\Char^\hLambda$ is a direct summand of $\CH(M)$ for some
irreducible $M\in\CM^\hLambda$, we see that $\HC(A)\ne0$, and hence $\hc(A)\ne0$. It follows that $\hc$ is
faithful. To prove that $\hc$ is essentially surjective and fully faithful we need to compare
$\CZ(\Har^\hlambda)$ with the center of a larger category of Harish-Chandra bimodules.

\subsection{Morita invariance of the center}\label{Mor_inv}
 We denote by $\Har^\hLambda$ the category of Harish-Chandra
bimodules over $U_\hLambda$ (finitely generated as a left $U$-module). 
Also, we denote by $\Har$ the category of (finitely
generated as a left $U$-module) Harish-Chandra bimodules (no restrictions on the action of $Z(U)$). Clearly, $\Har^\hLambda$
is a monoidal category containing
$U_\hlambda$ as a direct summand monoidal subcategory. Hence we have the functor of restriction to the
direct summand $Res:\ \CZ(\Har^\hLambda)\to\CZ(\Har^\hlambda)$.

\begin{lem} \label{morita} $Res:\ \CZ(\Har^\hLambda)\to\CZ(\Har^\hlambda)$ is an equivalence of braided
monoidal categories. \end{lem}

\proofpt 
Let $\Har^{\hatt{\la}, \hatt{\la+\La}}$ be the category of Harish-Chandra bimodules
finitely generated as $U(\g)$-modules,
where the left action of $U(\g)$ has generalized infinitesimal character $\la$.
(Notice that the right action then automatically factors through 
$\bigoplus\limits_{\mu\in \la+\La} U(\g)_{\hatt{\mu}}$).
The category  $\Har^{\hatt{\la}, \hatt{\la+\La}}$ carries an action of the tensor product
$\Har^{\hatt{\la}}\otimes \Har^{\hatt{\la+\La}}$, where the first factor acts by tensoring
over $U(\g)$ on the left, while the second one acts by tensoring over $U(\g)$ on the right.

We claim that the following ``double commutator" property holds.

{\em (DC) The functor from $\Har^{\hatt{\la+\La}}$ to the category of endo-functors
of $\Har^{\hatt{\la}, \hatt{\la+\La}}$ as an $\Har^{\hatt{\la}}$-module category
is a full embedding, its image consists of endo-functors killed by a finite codimension
ideal in the center of $U(\g)$.

Likewise, the functor from $\Har^{\hatt{\la}}$ to the category of endo-functors
of $\Har^{\hatt{\la}, \hatt{\la+\La}}$ as an $\Har^{\hatt{\la+\La}}$-module category
is a full embedding, its image consists of endo-functors killed by a power of the maximal
ideal $\la$ in the center of $U(\g)$.}

The double commutator property (DC) yields the Lemma by a standard argument
showing 2-Morita invariance of the categorical center (cf. \cite{O}):
it is clear from (DC) that both the center $\Z(\Har^{\hatt{\la}})$ and $\Z(\Har^{\hatt{\la+\La}})$
are equivalent to the category of $Z(U(\g))$-finite endo-functors of $\Har^{\hatt{\la}, \hatt{\la+\La}}$
as an $\Har^{\hatt{\la}}\otimes \Har^{\hatt{\la+\La}}$-module category.

It remains to prove (DC). 
Let $\End_{\Har^{\hatt{\la}}}^{\operatorname{fin}} ( \Har^{\hatt{\la}, \hatt{\la+\La}})$, 
$\End_{\Har^{\hatt{\la+\La}}}^{\operatorname{fin}} ( \Har^{\hatt{\la}, \hatt{\la+\La}})$ denote the category 
of $Z(U(\g))$ finite endo-functors of  $ \Har^{\hatt{\la}, \hatt{\la+\La}}$ as a
$\Har^{\hatt{\la}}$-module category (respectively, as a $ \Har^{\hatt{\la+\La}}$-module
category). We need to show that the natural functors
$$ \Har^{ \hatt{\la+\La}} \to \End_{\Har^{\hatt{\la}}}^{\operatorname{fin}} ( \Har^{\hatt{\la}, \hatt{\la+\La}})\eqno{(*)}$$
$$\Har^{\hatt{\la}}\to \End_{\Har^{\hatt{\la+\La}}}^{\operatorname{fin}} ( \Har^{\hatt{\la}, \hatt{\la+\La}})\eqno{(**)}$$
are equivalences.

We prove that the first functor is an equivalence by constructing the inverse 
equivalence. Let $F: \Har^{\hatt{\la},\hatt{\mu}}
\to \Har^{\hatt{\la},\hatt{\nu}}$ be a functor commuting with the
$\Har^{\hatt{\la}}$ action and annihilated by the $N$-th power of the maximal ideal $(\mu)$
in the center $Z(U(\g))$. 
We proceed to define an object $M\in \Har^{\hatt{\mu},\hatt{\nu}}$ with an isomorphism
between $F$ and the functor $X\mapsto X\otimes _{U} M$. 

Let $X\in  \Har^{\hatt{\la},\hatt{\mu}}$ be a Harish-Chandra bimodule, such that the right action
of $U(\g)$ factors through $U(\g)/(\mu)^N$ making it a projective generator for 
$U(\g)/(\mu)^N$-modules; for example, if $\la$ and $\mu$ are dominant (or, more generally,
if the sum of both $\la$ and  $\mu$ with any strictly dominant positive weight is regular)
 we can take $X$ to be the direct summand with the
specified action of the center 
in $V\otimes (U(\fg)/(\mu)^N)$ for a finite dimensional modules $V$ with an extremal
weight $\la-\mu$; a standard argument shows this module has the stated properties
as long as $\la$ is regular.  

The key step is to define a right action of the ring $\End_{U(\g)/(\mu)^N} (X)$ on the module $F(X)$
by a $U_{\hatt{\nu}}$- (though not by $U_{\hatt{\la}}$-) endomorphisms. Once this is done,
we can apply the Morita equivalence between $\End_{U(\g)/(\mu)^N} (X)^{op}$ and $U(\g)/(\mu)^N$
to $F(X)$ obtaining a $U(\g)/(\mu)^N$ module  $M$ with a commuting $U_{\hatt{\nu}}$-action;
we will then check that
the functor $F$ is given by $L\mapsto L\otimes _U M$, and also that
 the diagonal action of $\g$ on $M$ is locally finite, thus we get an object
 $M\in \Har^{\hatt{\mu},\hatt{\nu}}$. 
 
To construct the action, let us view the Harish-Chandra module $X$ as a 
$\g$-module
(where we use the right $\g$-action) equipped with a compatible $G$-action.
Then we have 
$$\Hom_\g(X,X) = \Hom_{\Har} (X, \CO(G)\otimes X) = 
\Hom_{\Har}(X, \bigoplus\limits_{\eta\in \Lambda^+} V_\eta
\otimes (V_\eta^*\otimes U(\g))\otimes_{U(\g)}X),$$
where $V_\eta$ denotes the irreducible $G$-module with highest weight $\eta$.

We can assume without loss of generality that $N$ is such that the left action of the 
ideal $(\la)^N$ in $Z(U)$ on $X$ vanishes; then we can rewrite
$\Hom_{\Har}(X,  (V_\eta^*\otimes U(\g))\otimes_{U(\g)}X)$ as 
$\Hom(X, M_\eta \otimes_{U(\g)}X)$,
where $M_\eta\in \Har^{\hatt{\la}}$ is the maximal summand in 
$(V_\eta^*\otimes U(\g))/(\la)^N$
where the left action of the center has generalized central character $\la$.
Since $F$ is compatible with the action of $\Har^{\hatt{\la}}$
we get the required map 
\begin{multline*}
\Hom_\g(X,X)= \bigoplus\limits_{\eta\in \Lambda^+}
\Hom_{\Har}(X,  V_\eta
\otimes (M_\eta \otimes_{U(\g)}X))\to\\ 
\to\bigoplus\limits_{\eta\in \Lambda^+}
\Hom_{\Har}(F(X),  V_\eta
\otimes (M_\eta \otimes_{U(\g)}F(X)) )
=\Hom_\g(F(X),F(X)).
\end{multline*}

Since $F(X)=X\otimes_U M$ for the  projective generator $X$, 
we see that $F$ is canonically isomorphic to the functor 
$L\mapsto L\otimes _U M$. Now Yoneda
Lemma implies that a bimodule $M$ with such an isomorphism is unique up to a unique
isomorphism if it exists (in particular it does not depend on the auxiliary choice
of the projective generator $X$).
Thus  the above procedure defines a functor from
$ \End_{\Har^{\hatt{\la}}}^{\operatorname{fin}} ( \Har^{\hatt{\la}, \hatt{\la+\La}})$
to $\g$-bimodules.

The diagonal action of $\g$ on the bimodule $M$ has to be locally finite:
to see this notice that $U(\g)/(\mu)^N$ is a quotient of a tensor product 
$P\otimes_{U(\g)}Q$ for some $P\in \Har^{\hatt{\mu},\hatt{\la}}$,
 $Q\in \Har^{\hatt{\la},\hatt{\mu}}$, thus $M$ is quotient of $P\otimes_{U(\g)}F(Q)$. 
 Thus we defined a functor $ \End_{\Har^{\hatt{\la}}}^{\operatorname{fin}} ( \Har^{\hatt{\la}, \hatt{\la+\La}})
 \to  \Har^{ \hatt{\la+\La}} $, it is easy to see that it is the equivalence inverse to $(*)$.

The construction of the equivalence inverse to $(**)$ is similar (though simpler).
 \qed

\subsection{Comparison to the center of $\Har$} \label{aq} Recall that $\Har$ 
is the category of (finitely
generated as a left $U$-module) Harish-Chandra bimodules 
(no restrictions on the action of $Z(U)$). Let
$(M,b_X)_{X\in\Har^\hLambda}$ be an object of $\CZ(\Har^\hLambda)$. We will construct a family of compatible
isomorphisms $b_Y:\ M\otimes_U Y\iso Y\otimes_U M$ for $Y\in\Har$, thus providing a fully faithful functor
$\CZ(\Har^\hLambda)\to\CZ(\Har)$.

To this end note that we have a canonical isomorphism $U_\hLambda\otimes_UY\simeq
Y\otimes_UU_\hLambda=:X(Y)\in\Har^\hLambda$. Now since $M\otimes_UY=M\otimes_{U_\hLambda}X(Y)$, and
$Y\otimes_UM=X(Y)\otimes_{U_\hLambda}M$, we can define $b_Y:=b_{X(Y)}$. It is immediate to extend the
assignment $(M,b_X)_{X\in\Har^\hLambda}\mapsto(M,b_Y)_{Y\in\Har}$ to a functor $\BF:\
\CZ(\Har^\hLambda)\to\CZ(\Har)$, and to prove that it is fully faithful (but certainly not essentially
surjective).

\subsection{The center of $\Har$ and $D$-modules on $G$} Recall that we have 2 commuting embeddings
$U\hookrightarrow\CalD(G)$ (as left invariant differential operators), $U^{opp}\simeq U\hookrightarrow
\CalD(G)$ (as right invariant differential operators). Thus the global sections of any $G$-equivariant
$\CalD(G)$-module $A$ have a structure of a Harish-Chandra bimodule, to be denoted by $\Gamma(G,A)$.
Moreover, $\Gamma(G,A)$ is equipped with a canonical central structure. Thus we obtain a functor
$\Gamma(G,-):\ \CalD(G)-mod^G\to\CZ(\Har)$.

\begin{lem} \label{benzvi} The functor $\Gamma(G,-):\ \CalD(G)-mod^G\to\CZ(\Har)$ is an equivalence of
abelian categories. \end{lem}

\proofpt Let $(M,b_Y)_{Y\in\Har}$ be an object of $\CZ(\Har)$. Let $Y=Fr(E)$ be a free Harish-Chandra
bimodule associated to a finite dimensional $G$-module $E$. We have $M\otimes_UY\simeq M\otimes E$, and
$Y\otimes_UM\simeq E\otimes M$. We will denote by $b_E$ the composition morphism $M\otimes E\simeq
M\otimes_UY\stackrel{b_Y}{\longrightarrow} Y\otimes_UM\simeq E\otimes M$ or, equivalently, $b_E:\ E\otimes
E^*\to \End(M)$. Note that for $E=E_1\otimes E_2$ the compatibility condition in the definition of Drinfeld
center implies that the composition $$\begin{CD} E_1\otimes E_2\otimes M @>{\Id_{E_1}\otimes b_{E_2}}>>
E_1\otimes M\otimes E_2 @>{b_{E_1}\otimes\Id_{E_2}}>> M\otimes E_1\otimes E_2 \end{CD}$$ equals $b_E$. In
other words, for any $v_{1,2}\in E_{1,2},\ v_{1,2}^*\in E_{1,2}^*$ we have an equality in the ring
$\End(M)$: $$b_E(v_1\otimes v_2\otimes v_1^*\otimes v_2^*)=b_{E_1}(v_1\otimes v_1^*) b_{E_2}(v_2\otimes
v_2^*).$$ Note that the same equality holds in the ring $\CO(G)$ for the matrix coefficients of the
representations $E,E_1,E_2$; moreover, $\CO(G)$ is generated by the matrix coefficients, and is given by the
above relations. Hence we obtain a homomorphism $\phi:\ \CO(G)\to\End(M)$. For $x\in\fg$ acting as a left
invariant vector field on $\CO(G)$ we have $\phi(xv\otimes v^*)(m)=x(\phi(v\otimes v^*)(m))$ for any $m\in
M,\ v\in E,\ v^*\in E^*$. In other words, the left $U$-action on $M$ together with the action of $\CO(G)$
combine into an action of $U\ltimes\CO(G)= \CalD(G)$. Now the action of $G$ on $M$ provides the resulting
$\CalD(G)$-module with a $G$-equivariant structure.

It is easy to see that this construction gives rise to a functor $F$ from $\CZ(\Har)$ to $\CalD(G)-mod^G$
quasiinverse to $\Gamma(G,-)$. \qed

\subsection{Completion of proof of Theorem~\ref{davids}} We consider the composition of functors
$$\begin{CD} \CZ(\Har^\hlambda) @>\bF>> \CZ(\Har^\hLambda) @>\BF>> \CZ(\Har) @>F>> \CalD(G)-mod^G \end{CD}$$
constructed respectively in the proof of Lemma~\ref{morita}, in subsection~\ref{aq}, and in the proof of
Lemma~\ref{benzvi}. It is easy to see that the composition lands into $\Char^\hLambda$, and is quasiinverse
to $\hc$. This completes the proof of~Theorem~\ref{davids}. \qed

\section{Truncated convolution and Harish-Chandra bimodules}

\subsection{Two-sided cells} For $w\in W$ let $\tO_w$ stand for the preimage in $\tB\times\tB$ of the
$G$-orbit ${\mathbb O}_w\subset\CB\times\CB$. Clearly, the restriction (both shriek and star) of any $M\in
D\CM^{\hLambda,\widehat{-\lambda+\Lambda}}$ to $\tO_w$ is zero if $w\not\in W_\zeta$ (the stabilizer of
$\zeta$ in $W$). Similarly, the restriction of any $M\in D\CM^{\hLambda,\widehat{-w_0\lambda+\Lambda}}$ to
$\tO_w$ is zero if $w\not\in w_0W_\zeta$. For $w\in W_\zeta$ we denote by $IC_w^{\zeta,\zeta^{-1}}$ the
isomorphism class of an irreducible $\CalD$-module in $\CM^{\hLambda,\widehat{-\lambda+\Lambda}}$ supported
on the closure of $\tO_w$.

Let $\uc\subset W_\zeta$ be a two-sided cell, and let $a(\uc)$ denote its $a$-function. Let
$\CM^{\hLambda,\widehat{-\lambda+\Lambda}}_{\leq\uc}$ stand for the Serre subcategory of
$\CM^{\hLambda,\widehat{-\lambda+\Lambda}}$ generated by the sheaves $IC_w^{\zeta,\zeta^{-1}}$ for $w$ lying
in $\uc$ and all the smaller cells. We denote by $\CM^{\hLambda,\widehat{-\lambda+\Lambda}}_\uc$ the
quotient of $\CM^{\hLambda,\widehat{-\lambda+\Lambda}}_{\leq\uc}$ by the Serre subcategory
$\CM^{\hLambda,\widehat{-\lambda+\Lambda}}_{<\uc}$ generated by the sheaves $IC_w^{\zeta,\zeta^{-1}}$ for
$w$ lying in all the cells smaller than $\uc$. We denote by
$D\CM^{\hLambda,\widehat{-\lambda+\Lambda}}_{\leq\uc}$ (resp.
$D\CM^{\hLambda,\widehat{-\lambda+\Lambda}}_{<\uc}$) the full triangulated subcategory of
$D\CM^{\hLambda,\widehat{-\lambda+\Lambda}}$ generated by
$\CM^{\hLambda,\widehat{-\lambda+\Lambda}}_{\leq\uc}$ (resp.
$\CM^{\hLambda,\widehat{-\lambda+\Lambda}}_{<\uc}$). We denote by
$D\CM^{\hLambda,\widehat{-\lambda+\Lambda}}_\uc$ the quotient of
$D\CM^{\hLambda,\widehat{-\lambda+\Lambda}}_{\leq\uc}$ by
$D\CM^{\hLambda,\widehat{-\lambda+\Lambda}}_{<\uc}$. It contains the heart
$\CM^{\hLambda,\widehat{-\lambda+\Lambda}}_\uc$.

We introduce the similar categories and notations for $D\CM^{\hLambda,\widehat{-w_0\lambda+\Lambda}}$ in
place of $D\CM^{\hLambda,\widehat{-\lambda+\Lambda}}$. The only difference is in the numeration of
irreducible objects: now $IC_w^{\zeta,w_0\zeta^{-1}}$ stands for the isomorphism class of an irreducible
$\CalD$-module in $\CM^{\hLambda,\widehat{-w_0\lambda+\Lambda}}$ supported on the closure of $\tO_{ww_\zeta
w_0}$. Here $w_\zeta$ stands for the longest (in $W$) element of $W_\zeta$, so that $w_\zeta w_0$ is the
shortest element in the left coset $W_\zeta w_0\subset W$. With these notations the intertwining functor
$\bI_{w_0}$ respects the filtrations by the cell subcategories:
$\bI_{w_0}(D\CM^{\hLambda,\widehat{-\lambda+\Lambda}}_{\leq\uc})\subset
D\CM^{\hLambda,\widehat{-w_0\lambda+\Lambda}}_{\leq\uc}$, and gives rise to the same named functor
$\bI_{w_0}:\ D\CM^{\hLambda,\widehat{-\lambda+\Lambda}}_\uc\to
D\CM^{\hLambda,\widehat{-w_0\lambda+\Lambda}}_\uc$.

For $w\in W_\zeta$, we denote by $V^\lambda_w$ the isomorphism class of an irreducible Harish-Chandra
bimodule such that $V^\lambda_w\simeq R\Gamma^{\hlambda,\widehat{-w_0\lambda}}(IC_w^{\zeta,w_0\zeta^{-1}})$.
Now for a two-sided cell $\uc\in W_\zeta$ we define the categories
$\Har^\hlambda_{<\uc}\subset\Har^\hlambda_{\leq\uc} \subset\Har^\hlambda,\ \Har^\hlambda_\uc,\
D\Har^\hlambda_\uc$, etc. in an evident fashion. The global sections functors
$R\Gamma^{\hlambda,\widehat{-w_0\lambda}}:\ D\CM^{\hLambda,\widehat{-w_0\lambda+\Lambda}}\to
D\Har^\hlambda;\ R\Gamma^{\hlambda,\hulambda}:\ D\CM^{\hLambda,\widehat{-\lambda+\Lambda}}\to
D\Har^\hlambda$ preserve the filtrations by the cell subcategories, and give rise to the same named functors
$R\Gamma^{\hlambda,\widehat{-w_0\lambda}}:\ D\CM^{\hLambda,\widehat{-w_0\lambda+\Lambda}}_\uc\to
D\Har^\hlambda_\uc;\ R\Gamma^{\hlambda,\hulambda}:\ D\CM^{\hLambda,\widehat{-\lambda+\Lambda}}_\uc\to
D\Har^\hlambda_\uc$.

\begin{prop} \label{mit} The functor $\bI_{w_0}[a_\uc]:\ D\CM^{\hLambda,\widehat{-\lambda+\Lambda}}_\uc\to
D\CM^{\hLambda,\widehat{-w_0\lambda+\Lambda}}_\uc$ is exact, that is
$\bI_{w_0}[a_\uc](\CM^{\hLambda,\widehat{-\lambda+\Lambda}}_\uc)
\subset\CM^{\hLambda,\widehat{-w_0\lambda+\Lambda}}_\uc$. \end{prop}

\begin{rem} \label{512}
 The Proposition can be viewed as a categorification of Lusztig's
result~\cite[5.12.2]{L84}. The latter shows that the action of the element $T^2_{w_0}$ on the two-sided cell
subquotient module $H_\uc$ of the Hecke algebra $H$ is the scalar multiplication by $q^{d/2}$ for some
integer $d=d(\uc)$. More precisely,
 the Proposition implies that the functor $\bI_{w_0}^{-2}[-2 a_\uc]$ induces
an exact functor on the cell subquotient. Since the functor lifts to the graded version of the category and
the induced map on the (non-graded) Grothendieck group is trivial, it follows that the map induced by
$T_{w_0}^2$ on $H_\uc$ sends an element $C_w$ of the
  Kazhdan-Lusztig basis to $q^{d(w)/2}C_w$. In fact,
in {\em loc. cit.} Lusztig computes the value of $d(w)=d(\uc)$. It would be interesting to understand the
answer from the point of view of the present paper. Notice that $d(w)$ does {\em not} coincide with
$-a_\uc$, which agrees with the fact that the functor $\bI_{w_0}^{-1}$ is not pure.
\end{rem}

\begin{rem}
By the result of \cite{BBM}, the functor $\bI_{w_0}^2$ is the {\em Serre functor}
for the derived category $D_N(G/B)$. It can be deduced from the Proposition that
the functor on the cell subquotient induced by $\bI_{w_0}^2$ is isomorphic to the functor of shift by $-2a_{\uc}$, i.e. the cell subquotient  is a Calabi--Yau category of dimension $2a_{\uc}$. It is interesting to compare it with the results
of Mazorchuk and Stroppel \cite{Str1}, \cite{Str2} which imply that for groups of type $A_n$ the (right) cell
subquotients in the {\em abelian} category of perverse sheaves are equivalent to the category of modules over a  symmetric (in an alternative terminology, Frobenius) algebra, thus are Calabi-Yau of dimension zero.
\end{rem}

\begin{rem}\label{mir}
The Calabi--Yau property of the cell subquotient categories was conjectured by Ivan Mirkovi\' c (and communicated to one of us
around 2005).
\end{rem}

\begin{cor} \label{trunc} a) The functor $R\Gamma^{\hlambda,\hulambda}[a_\uc]:\
D\CM^{\hLambda,\widehat{-\lambda+\Lambda}}_\uc\to D\Har^\hlambda_\uc$ is exact, that is
$R\Gamma^{\hlambda,\hulambda}[a_\uc] (\CM^{\hLambda,\widehat{-\lambda+\Lambda}}_\uc)\subset
\Har^\hlambda_\uc$.

b) This functor is monoidal with respect to the truncated convolution on
$\CM^{\hLambda,\widehat{-\lambda+\Lambda}}_\uc:\ (M_1,M_2)\mapsto M_1\bullet M_2:=\ul{H}^{a_\uc}(M_1*M_2)$,
and the nonderived convolution on $\Har^\hlambda_\uc:\ (V_1,V_2)\mapsto V_1\otimes_U V_2$. \end{cor}

\proofpt a) is immediate, while b) follows from a) and the monoidal property of $
R\Gamma^{\hlambda,\hulambda}:\ D\CM^{\hLambda,\widehat{-\lambda+\Lambda}}\to D\Har^\hlambda$. \qed

The proof of the proposition occupies the rest of the section.

\begin{rem}\label{proJ}
 Notice that the last claim
 shows that Lusztig's asymptotic Hecke algebra $J$ can be interpreted
as the Grothendieck group of the monoidal category of semi-simple
 Harish-Chandra bimodules
where the monoidal structure comes from tensor product of bimodules
 taken modulo
the subcategory generated by the submodules with a smaller support.
A related  statement was obtained by A.~Joseph \cite{J1}.
\end{rem}

\subsection{Category $\CO$} \label{O} We denote by
$\bM^{\widehat{-\lambda+\Lambda},\widehat{-\lambda+\Lambda}}$ the category of $N$-equivariant,
weakly
$B$-equivariant $\CalD_\tB$-modules such that the locally finite action of $\BC[\ft^*]$ (differential
operators arising from the infinitesimal action of $T\subset B\subset G$ on $\tB$) has generalized
eigenvalues in $-\lambda+\Lambda$; while the action of $\BC[\ft^*]$ (differential operators arising from the
infinitesimal right action of $T$ on $\tB$) is locally finite with generalized eigenvalues in
$-\lambda+\Lambda$. For $w\in W_\zeta$, we denote by $IC^{\zeta^{-1}}_w$ the isomorphism class of an
irreducible $\CalD_\tB$-module in $\bM^{\widehat{-\lambda+\Lambda},\widehat{-\lambda+\Lambda}}$ supported on
the closure of $\tB_w$ (preimage in $\tB$ of the corresponding $B$-orbit in $\CB$). We denote by
$D\bM^{\widehat{-\lambda+\Lambda},\widehat{-\lambda+\Lambda}}$ the derived category of
$\bM^{\widehat{-\lambda+\Lambda},\widehat{-\lambda+\Lambda}}$. The categories
$\bM^{\widehat{-\lambda+\Lambda},\widehat{-w_0\lambda+\Lambda}}\subset
D\bM^{\widehat{-\lambda+\Lambda},\widehat{-w_0\lambda+\Lambda}}$ are defined similarly: we require that the
generalized eigenvalues of the right action of $\BC[\ft^*]$ lie in $-w_0\lambda+\Lambda$. For $w\in
W_\zeta$, we denote by $IC^{w_0\zeta^{-1}}_w$ the isomorphism class of an irreducible $\CalD_\tB$-module in
$\bM^{\widehat{-\lambda+\Lambda},\widehat{-w_0\lambda+\Lambda}}$ supported on the closure of $\tB_{ww_\zeta
w_0}$.

We denote by $\CO^\hulambda$ the category of $N$-integrable $U$-modules where $Z(U)$ acts locally finitely
with the generalized eigenvalue $-\lambda-2\rho$ (same as $-w_0\lambda$), and $\ft$ acts locally finitely
with the generalized eigenvalues in $-\lambda+\Lambda$. We denote by $D\CO^\hulambda$ the derived category
of $\CO^\hulambda$.

For $M\in\bM^{\widehat{-\lambda+\Lambda},\widehat{-\lambda+\Lambda}}$ we denote by $\Gamma^\hulambda(\tB,M)$
the maximal direct summand of the global sections on which $\BC[\ft^*]$ acts with the generalized eigenvalue
$-\lambda-2\rho$. It gives rise to an equivalence of categories $R\Gamma^\hulambda(\tB,-):\
D\bM^{\widehat{-\lambda+\Lambda}, \widehat{-\lambda+\Lambda}}\to D\CO^\hulambda$. Similarly, we have an
equivalence of categories $R\Gamma^{\widehat{-w_0\lambda}}(\tB,-):\
D\bM^{\widehat{-\lambda+\Lambda},\widehat{-w_0\lambda+\Lambda}}\to D\CO^\hulambda$. The latter equivalence
is exact, that is, takes $\bM^{\widehat{-\lambda+\Lambda},\widehat{-w_0\lambda+\Lambda}}$ to
$\CO^\hulambda$.

We have a closed embedding $\iota:\ \tB\to\tB\times\tB,\ x\mapsto(e,x)$ where $e\in\tB$ is the image of the
neutral element $e\in G$ under $G\to G/N=\tB$. The functor $\iota^o:=\iota^![\dim\tB]$ gives rise to the
same named equivalences $\CM^{\hLambda,\widehat{-\lambda+\Lambda}}\iso
\bM^{\widehat{-\lambda+\Lambda},\widehat{-\lambda+\Lambda}};\
\CM^{\hLambda,\widehat{-w_0\lambda+\Lambda}}\iso
\bM^{\widehat{-\lambda+\Lambda},\widehat{-w_0\lambda+\Lambda}}$. We keep the same name $\iota^o$ for the
derived versions of the above functors. All the notations related to two-sided cells are transfered to the
categories $\bM$ via the functors $\iota^o$.

For $w\in W_\zeta$, we denote by $L_w^{-\lambda-2\rho}$ the isomorphism class of an irreducible module in
$\CO^\hulambda$ with the highest weight $w\cdot(-\lambda-2\rho)$. Note that
$R\Gamma^{\widehat{-w_0\lambda}}(\tB,IC_w^{w_0\zeta^{-1}}) \simeq L_w^{-\lambda-2\rho}$. Now for a two-sided
cell $\uc\subset W_\zeta$ we define the categories $\CO^\hulambda_{<\uc}\subset\CO^\hulambda_{\leq\uc}
\subset\CO^\hulambda,\ \CO^\hulambda_\uc,\ D\CO^\hulambda_\uc$, etc. in an evident fashion.

We have a natural isomorphism of functors $\bI_{w_0}^{-1}\circ\iota^o\simeq\iota^o\circ\bI_{w_0}^{-1}:\
D\CM^{\hLambda,\widehat{-w_0\lambda+\Lambda}}\to
D\bM^{\widehat{-\lambda+\Lambda},\widehat{-\lambda+\Lambda}}$. Here $\bI_{w_0}^{-1}:\
D\bM^{\widehat{-\lambda+\Lambda},\widehat{-w_0\lambda+\Lambda}}\to
D\bM^{\widehat{-\lambda+\Lambda},\widehat{-\lambda+\Lambda}}$ is
as in section \ref{inter}.

Finally, we have a natural isomorphism of functors $$ R\Gamma^\hulambda(\tB,-)\circ\bI_{w_0}^{-1}\simeq
R\Gamma^{\widehat{-w_0\lambda}}(\tB,-):\ D\bM^{\widehat{-\lambda+\Lambda},\widehat{-w_0\lambda+\Lambda}}\to
D\CO^\hulambda$$ (cf.~\cite{BB}). It follows that the functor $R\Gamma^\hulambda(\tB,-)$ respects the
filtrations by the cell subcategories: $R\Gamma^\hulambda(\tB,-)(D\bM^{\widehat{-\lambda+\Lambda},
\widehat{-\lambda+\Lambda}}_{\leq\uc}) \subset D\CO^\hulambda_{\leq\uc}$, and gives rise to the same named
functor $R\Gamma^\hulambda(\tB,-):\ D\bM^{\widehat{-\lambda+\Lambda},\widehat{-\lambda+\Lambda}}_\uc\to
D\CO^\hulambda_\uc$.

Clearly, Proposition~\ref{mit} is equivalent to the following

\begin{prop} \label{mit'} The functor $\bI_{w_0}[a_\uc]:\
D\bM^{\widehat{-\lambda+\Lambda},\widehat{-\lambda+\Lambda}}_\uc\to
D\bM^{\widehat{-\lambda+\Lambda},\widehat{-w_0\lambda+\Lambda}}_\uc$ is exact, that is
$\bI_{w_0}[a_\uc](\bM^{\widehat{-\lambda+\Lambda}, \widehat{-\lambda+\Lambda}}_\uc)
\subset\bM^{\widehat{-\lambda+\Lambda},\widehat{-w_0\lambda+\Lambda}}_\uc$. Equivalently, the functor
$R\Gamma^\hulambda(\tB,-)[a_\uc]:\ D\bM^{\widehat{-\lambda+\Lambda},\widehat{-\lambda+\Lambda}}_\uc\to
D\CO^\hulambda_\uc$ is exact, that is $R\Gamma^\hulambda(\tB,-)[a_\uc]
(\bM^{\widehat{-\lambda+\Lambda},\widehat{-\lambda+\Lambda}}_\uc)\subset \CO^\hulambda_\uc$. \end{prop}

The rest of the section is devoted to the proof.

\subsection{Dualities} \label{dual} For a $U$-module $M$, there is a natural right $U$-module structure on
$\Hom_U(M,U)$. If $M\in\CO^\hulambda$, the center $Z(U)$ acts on (the cohomology modules of) $\RHom_U(M,U)$
locally finitely with generalized eigenvalues $-\lambda-2\rho$. Composing with the opposition isomorphism
$U\simeq U^{opp}$, we get a duality functor $\RHom_U^{opp}(-,U):\ D\CO^\hulambda\to (D\CO^\hlambda)^{opp}$.
The same construction gives rise to the quasiinverse duality functor $\RHom_U^{opp}(-,U):\ D\CO^\hlambda\to
(D\CO^\hulambda)^{opp}$. For $w\in W_\zeta$, we denote by $L_w^\lambda$ the isomorphism class of an
irreducible module in $\CO^\hlambda$ with the highest weight $ww_\zeta\cdot\lambda$. For a two-sided cell
$\uc\in W_\zeta$ we define the categories $\CO^\hlambda_{<\uc}\subset\CO^\hlambda_{\leq\uc}
\subset\CO^\hlambda,\ \CO^\hlambda_\uc,\ D\CO^\hlambda_\uc$, etc. in an evident fashion. The duality
functors preserve the filtrations by the cell subcategories and give rise to the same named functor
$\RHom_U^{opp}(-,U):\ D\CO^\hlambda_\uc\to (D\CO^\hulambda_\uc)^{opp}$.

For a $\CalD_\tB$-module $M$, there is a natural right $\CalD_\tB$-module structure on ${\mathcal
H}om_{\CalD_\tB}(M,\CalD_\tB)$. If $M\in\bM^{\widehat{-\lambda+\Lambda},\widehat{-\lambda+\Lambda}}$, the
action of $\BC[\ft^*]$ (coming from the infinitesimal right action of $T$ on $\tB$) on (the cohomology
sheaves of) $R{\mathcal H}om_{\CalD_\tB}(M,\CalD_\tB)$ is locally finite with generalized eigenvalues in
$-\lambda+\Lambda$. We have a canonical opposition isomorphism $\CalD_{\omega_\tB}\simeq\CalD_\tB^{opp}$
where $\omega_\tB$ stands for the canonical line bundle on $\tB$, and $\CalD_{\omega_\tB}$ denotes the
corresponding tdo. We have an isomorphism $\omega_\tB\simeq\CO_\tB$ which gives rise to the isomorphism
$\CalD_{\omega_\tB}\simeq\CalD_\tB$; however the natural action of $\BC[\ft^*]$ gets shifted by $2\rho$
under this isomorphism. Thus we obtain a duality functor $\BD:=R{\mathcal H}om_{\CalD_\tB}(-,\CalD_\tB):\
D\bM^{\widehat{-\lambda+\Lambda},\widehat{-\lambda+\Lambda}}\to (D\bM^{\hLambda,\hLambda})^{opp}$.

By construction, we have a natural isomorphism of functors $R\Gamma^\hlambda(\tB,-)\circ\BD\simeq
\RHom_U^{opp}(-,U)\circ R\Gamma^\hulambda:\ D\bM^{\widehat{-\lambda+\Lambda},\widehat{-\lambda+\Lambda}}\to
(D\CO^\hlambda)^{opp}$.

Since any $M\in\bM^{\widehat{-\lambda+\Lambda},\widehat{-\lambda+\Lambda}}$ is holonomic, $\BD[\dim\tB]$ is
exact, that is, takes $\bM^{\widehat{-\lambda+\Lambda},\widehat{-\lambda+\Lambda}}$ to
$(\bM^{\hLambda,\hLambda})^{opp}$. In particular, for $w\in W_\zeta$, we have $\BD
IC_w^{\zeta^{-1}}[\dim\tB]\simeq IC_w^\zeta$ where $IC_w^\zeta$ is the isomorphism class of an irreducible
$\CalD_\tB$-module in $\bM^{\hLambda,\hLambda}$ supported on the closure of $\tB_w$. Moreover,
$R\Gamma^\hlambda(\tB,-):\ D\bM^{\hLambda,\hLambda}\to D\CO^\hlambda$ is exact, that is, takes
$\bM^{\hLambda,\hLambda}$ to $\CO^\hlambda$. In particular, $R\Gamma^\hlambda(\tB,IC_w^\zeta)\simeq
L^\lambda_w$. Hence the next lemma implies that the functor $R\Gamma^\hulambda(\tB,-)[a_\uc]:\
D\bM^{\widehat{-\lambda+\Lambda},\widehat{-\lambda+\Lambda}}_\uc\to D\CO^\hulambda_\uc$ is  
exact, thus it yields Proposition \ref{mit'}.

\begin{lem} \label{ias} The functor $\RHom_U^{opp}(-,U)[\dim\tB+a_\uc]:\ D\CO^\hlambda_\uc\to
(D\CO^\hulambda_\uc)^{opp}$ is 
 exact.
\end{lem}

\proofpt Recall the notion of the {\em associated variety} (see e.g.~\cite{J}). For $M\in\CO^\hlambda$ the
associated variety $V(M)$ is a closed subvariety of $\fn\subset\fg=\fg^*$ (here $\fn$ is the Lie algebra of
$N$, and we have identified $\fg$ with $\fg^*$ via the Killing form). If $w\in W_\zeta$, and
$M=L_w^\lambda$, then according to~3.10,~3.11 of {\em loc. cit.}, there is a nilpotent orbit
$\CN_{w,\lambda}\subset\CN\subset\fg$ with the closure $\overline\CN_{w,\lambda}$ such that any irreducible component of $V(L_w^\lambda)$ is an irreducible component of $\fn\cap\overline\CN_{w,\lambda}$. Moreover, if $w,w'\in\uc$, then $\CN_{w,\lambda}=\CN_{w',\lambda}=: \CN_\uc$, and $\dim\CN-\dim\CN_\uc=2m_\uc$, where $m_\uc$ is the degree of the Goldie rank polynomial associated to $\on{Ann}(L_w^\lambda)$ for any $w\in\uc$
(see~\cite{Jo},~Theorem~5.1). Finally, it is well known, due to the works of A.~Joseph and G.~Lusztig, that
$m_\uc=a_\uc$ in case of integral $\lambda\in\Lambda$. The equality $m_\uc=a_\uc$ in case of general
$\lambda$ is reduced to the case of integral $\lambda$ in~\cite{ABV} (Chapters 16, 17).

Since $\fn\cap\CN_\uc$ is Lagrangian in $\CN_\uc$ we see that for $w\in\uc$, the associated variety
$V(L_w^\lambda)$ is equidimensional of codimension $a_\uc$ in $\fn$, that is of codimension $a_\uc+\dim\tB$
in $\fg$. It follows (see e.g.~Theorem~V.2.2.2 of~\cite{B}) that for $i\ne 0$ the module
$\RHom_U^i(L_w^\lambda,U)^{opp}[\dim\tB+a_\uc]$ has associated variety of a smaller dimension.
 The lemma follows. \qed

\section{Classification and convolution of character $D$-modules} 

\subsection{Cells in character sheaves}\label{cCS} Recall that $\lambda$ is a regular weight.

For $w\in W_\zeta$, the irreducible objects $IC^{\zeta,\zeta^{-1}}_w$ of
$\CM^{\hLambda,\widehat{-\lambda+\Lambda}}$ have a natural equivariant structure with respect to the
diagonal right action of $T$ on $\tB\times\tB$, so by abuse of notation we will consider
$IC^{\zeta,\zeta^{-1}}_w$ as an isomorphism class of irreducible objects of $\CM^\hLambda$. For a two-sided
cell $\uc\in W_\zeta$ we define the categories $\CM^\hLambda_{<\uc}\subset\CM^\hLambda_{\leq\uc}
\subset\CM^\hLambda,\ \CM^\hLambda_\uc,\ D\CM^\hLambda_\uc$, etc. in an evident fashion.

We let $\Char^\hLambda_{<\uc}$ (respectively, $\Char^\hLambda_{\leq\uc}$) be the Serre subcategory in
$\Char^\hLambda$ given by $M\in \Char^\hLambda_{<\uc}$ (respectively, $M\in\Char^\hLambda_{\leq\uc}$)
 if $\Gamma^\hlambda(M) \in \Har^\hlambda_{<\uc}$
(respectively,  $\Har^\hlambda_{\leq \uc}$); see section~\ref{huji} for the definition of $\Gamma^\hlambda$.
It is easy to show that the subcategories do not depend on the choice of a regular $\lambda$ in a given coset
of $\Lambda$.

We define $D\Char^\hLambda_{<\uc}\subset D\Char^\hLambda$ as the full
 triangulated subcategory generated by $\Char^\hLambda_{<\uc}$, and similarly for $\leq\uc$; these are full
triangulated subcategories in  $D\Char^\hLambda$.
In view of the tensor property of the functor $\Gamma^\hlambda$ (which follows from Lemma \ref{conv_abs}
applied to the action of $G$ on itself by left translations) $D\Char^\hLambda_{<\uc}$,
$D\Char^\hLambda_{\leq \uc}$ are tensor ideals in $D\Char^\hLambda$.

 We denote by $\Char^\hLambda_\uc$ (respectively,
$D\Char^\hLambda_\uc$) the quotient category $\Char^\hLambda_{\leq\uc}/\Char^\hLambda_{<\uc}$ (respectively,
$D\Char^\hLambda_{\leq\uc}/D\Char^\hLambda_{<\uc}$). The functors $\HC$ and $\CH$ give rise to the same
named functors $\HC:\ D\Char^\hLambda_\uc\to D\CM^\hLambda_\uc;\ \CH:\ D\CM^\hLambda_\uc\to
D\Char^\hLambda_\uc$.

Recall the formula~(\ref{HoKa}). For $V\in\Har^\hlambda$ we define $L\Ind(V):=\CalD_\hLambda(G)\stackrel
{L}{\otimes}_{U_\hlambda\otimes U_\hulambda}~V$. This gives rise to a functor $L\Ind:\ D\Har^\hlambda\to
D\Char^\hLambda$ with a natural isomorphism of functors $L\Ind\circ R\Gamma^{\hlambda,\hulambda}\simeq\CH:\
D\CM^\hLambda\to D\Char^\hLambda$.

We have a natural isomorphism of functors
$R\Gamma^{\hlambda,\hulambda}\circ\on{p}^o\circ\HC\simeq\Gamma^\hlambda(G,-):\ D\Char^\hLambda\to
D\Har^\hlambda$. The functors $L\Ind$ and $\Gamma^\hlambda(G,-)$ preserve the filtrations by the cell
subcategories and give rise to the same named functors $L\Ind:\ D\Har^\hlambda_\uc\to D\Char^\hLambda_\uc;\
\Gamma^\hlambda(G,-):\ D\Char^\hLambda_\uc\to D\Har^\hlambda_\uc$. The functor $L\Ind$ is clearly right
exact; we denote by $\Ind$
its 0th cohomology sheaf.
The standard adjointness between the functors of restriction and induction shows that
the functor $\Ind:\ \Har^\hlambda\to \Char^\hLambda$ is left adjoint to the functor
$\Gamma^\hlambda(G,-)$ and same for their derived functors. Furthermore, it is easy to see that given a pair of adjoint triangulated
 functors
$F,\,G$
between triangulated
 categories $\A$, $\B$ and thick subcategories $\A'\subset \A$, $\B'\subset \B$
 such that $F:\A'\to \B'$, $G:\B'\to \A'$
 we get adjoint functors between $\A/\A'$ and $\B/\B'$. Thus we 
 see that  $L\Ind:\ D\Har^\hlambda_\uc\to D\Char^\hLambda_\uc$ is left adjoint to
$\Gamma^\hlambda(G,-):\ D\Char^\hLambda_\uc\to D\Har^\hlambda_\uc$,
hence $\Ind: \Har^\hlambda_\uc\to \Char^\hLambda_\uc$ is left adjoint to
$\Gamma^\hlambda(G,-):\ \Char^\hLambda_\uc\to \Har^\hlambda_\uc$.


We summarize the results of previous sections in the following

\begin{prop} \label{exact} a) The functor $\Gamma^\hlambda(G,-):\ D\Char^\hLambda_\uc\to D\Har^\hlambda_\uc$
is exact, that is $\Gamma^\hlambda(G,-)(\Char^\hLambda_\uc)\subset\Har^\hlambda_\uc$.

b) The functor $\CH[a_\uc]:\ D\CM^\hLambda_\uc\to D\Char^\hLambda_\uc$
sends $D^{\leq 0}\CM^\hLambda_\uc$ to $D^{\leq 0}\Char^\hLambda_\uc$, and for
$M\in\CM^\hLambda$, the truncation $\ul{H}^{a_\uc}\CH(M)$ is canonically isomorphic to
$\Ind(R^{a_\uc}\Gamma^{\hlambda,\hulambda}(M))\mod  \Char^\hLambda_{<\uc}$. \end{prop}

\proofpt a) is contained in Proposition~\ref{HotKas}.

Isomorphism \eqref{HoKa} shows that $\CH(M)\cong LInd(R\Gamma^{\hlambda,\hulambda}(M))$
for all $M\in \CM^\hLambda$.
Assuming that $M\in \CM^\hLambda_{\leq \uc}$, the image of
$R\Gamma^{\hlambda,\hulambda}(M)$ in the quotient category by smaller cells
is concentrated in homological degree
$a_\uc$ by Corollary~\ref{trunc}.a. Thus statement (b) follows by right exactness
of the induction functor.  \qed

As follows from the discussion in \ref{Drcsec} compared with Proposition \ref{HotKas}(a),
the functor $\Gamma^\hlambda(G,-):\ \Char^\hLambda_\uc\to\Har^\hlambda_\uc$ has a canonical central
structure and gives rise to the functor $\hc:\ \Char^\hLambda_\uc\to\CZ(\Har^\hlambda_\uc)$.


\subsection{Semisimplified categories}
We denote by $\Char^\hLambda_\ucs$ the full subcategory of
$\Char^\hLambda_\uc$ formed by all the semisimple objects. We denote by $\Har^\hlambda_\ucs$ (resp.
$\CM^\hLambda_\ucs$) the full subcategory of $\Har^\hlambda_\uc$ (resp. $\CM^\hLambda_\uc$) formed by all
the semisimple objects.
The functor $\CH$ is defined as a pull-back under a smooth map followed by push-forward
under a proper map (see section \ref{CSdef}), thus Decomposition Theorem \cite{BBD}
shows that it sends semi-simple objects to semisimple complexes. In particular,
 for $M\in\CM^\hLambda_\ucs$, the truncation
$\ul{H}^{a_\uc}\CH(M)$ lies in $\Char^\hLambda_\ucs$.
 Hence Proposition~\ref{exact} implies that for
$V\in\Har^\hlambda_\ucs$ the truncation $\Ind(V)$ lies in $\Char^\hLambda_\ucs$. We obtain the 
functor $\Ind_\uc:\ \Har^\hlambda_\ucs\to\Char^\hLambda_\ucs$. Recall that $\CM^\hLambda$ is identified with a
full subcategory of $\CM^{\hLambda,\widehat{-\lambda+\Lambda}}$. It follows that
$\CM^\hLambda_\ucs\simeq\CM^{\hLambda,\widehat{-\lambda+\Lambda}}_\ucs$. Clearly, all the irreducibles in
$\CM^\hLambda \simeq\CP^\hzeta$ lie in $\CP^\zeta$.
Since the convolution map for the flag variety is proper,
 the convolution of two
simple objects in $\CP^\zeta$ is a direct sum of shifts of simple objects
by Decomposition Theorem. Hence the {\em top cohomology
sheaf} of the convolution of two simple objects of $\CM^{\hLambda,\widehat{-\lambda+\Lambda}}$ is
semisimple. It follows that the category
$\CM^\hLambda_\ucs\simeq\CM^{\hLambda,\widehat{-\lambda+\Lambda}}_\ucs$ is monoidal with respect to the
operation $(M_1,M_2)\mapsto M_1\bullet M_2$ (see Corollary~\ref{trunc}.b); equivalently,
$\Har^\hlambda_\ucs$ is monoidal with respect to the operation $(V_1,V_2)\mapsto V_1\otimes_UV_2$. Both
categories $\CM^\hLambda_\ucs$ and $\Har^\hlambda_\ucs$ are equivalent to the semisimple monoidal category
${\mathcal C}_\zeta^\uc$ of~\cite{BFO}~\S5.

Recall that a  {\em commutator structure} on a functor $F$ from a monoidal category $\A$ to a category $\B$
is an isomorphism $F(M N)
\cong F (NM)$ fixed for all $M,N\in \A$  and
 satisfying the natural compatibilities, see~\cite{BFO}\S6 for details.

\begin{lem}\label{comstr}
The functors $\Ind:D^b(\Har^\hlambda) \to D^b( \Char^\hLambda_\ucs)$,
 $\Ind_\uc: \Har^\hlambda_\ucs\to \Char^\hLambda_\ucs$ carry a natural commutator structure.
\end{lem}

\proofpt By Theorem \ref{davids} the functor $\hc$ is a central functor. Hence the adjoint functor
$\Ind$ carries a commutator structure. In \cite[Proposition 5]{BFO} this implication is stated under
 stronger assumptions on the monoidal categories; however, it is also easy to show provided the target
 category of the central functor has a weak rigidity. The category of Harish-Chandra bimodules has a weak
 rigidity which sends a Harish-Chandra bimodule $M$ to the sum of  dual spaces to the isotypic
components of the action of $G$ on $M$ (notice that if $M$ is finitely generated and has a (generalized)
central character then the isotypic components are finite dimensional.)

 The commutator structure on $\Ind$ induces one on $\Ind_\uc$. \qed

Recall the following general statement. Let $\A$ be a rigid
 semisimple monoidal category with a finite number of isomorphism
classes of irreducible objects over an algebraically closed field of characteristic zero.
 Let $F: \CZ(\A)\to \A$ be the forgetful
functor and let $I: \A \to \CZ(\A)$ be its right adjoint (it exists since by \cite[Theorem 2.15]{ENO} the
category $\CZ(\A)$ is semisimple with a finite number of isomorphism classes of irreducible objects).
Then  the
functor $I$ has a natural structure of a commutator functor and is universal among commutator functors, see
\cite[\S6]{BFO} and also Lemma \ref{commz} below.

Thus the commutator functor $\Ind_\uc$ gives rise to a functor $\fA_\uc:
\CZ(\Har^\hlambda_\ucs)\to\Char^\hLambda_\ucs$.

\begin{thm} \label{december} a) The functor $\fA_\uc:\ \CZ(\Har^\hlambda_\ucs)\to\Char^\hLambda_\ucs$ is an
equivalence of categories.

b) The functor $\fA_\uc$ has a natural tensor structure,
 where $\Char^\hLambda_\ucs$ is equipped with the monoidal
  structure via the truncated convolution $A\circ B:=\underline{H}^0(A*B)$.
 Thus $\fA_\uc$ is a tensor equivalence. \end{thm}

 The Theorem will be proved in subsection~\ref{54}.

The Theorem together with the results of \cite{BFO} allows one to reprove the classification of irreducible
character sheaves (in the particular case of characteristic zero ground field) obtained by
Lusztig~\cite{L85} in a totally different way, and also to get new information about convolution of character sheaves. To spell this out recall that the structure of truncated
convolution categories was described in~\cite{BFO}\footnote{A more conceptual approach
to this description will be presented in a forthcoming paper \cite{LO}. There the finite set involved in Lusztig's conjectural description is interpreted as the set of finite dimensional modules over the finite $W$-algebra with a fixed central character, and the description of the truncated convolution category is derived from the action of the category of Harish-Chandra bimodules on the category of representations of the $W$-algebra.}
 following conjectures of Lusztig under the following
technical assumption (see~Conjecture~1 and~Theorems~3,4 of {\em loc. cit.}).

{\em (!) Either $G$ is of type $A$, or the centralizer of monodromy $\zeta$ in the Langlands dual group
$\check G$ is connected (which is automatically true for the adjoint $G$).}

For each two-sided cell $\uc$ Lusztig \cite{L84} defined a certain finite group ${\mathcal G}_{\uc}$ (this group is a quotient of the group of components of the centralizer of the corresponding special element).

Let us say that an equivalence $F$ of braided categories $\CC_1$, $\CC_2$ is a quasitensor equivalence if for every $M,N\in Ob(\CC_1)$ there exists an isomorphism
$F(MN)\cong F(M)F(N)$ (the existence of a choice of the isomorphisms satisfying the natural compatibilities is not required).

\begin{cor}\label{cor54} a) Under the assumption (!), the set of isomorphism classes of irreducible objects in
$\Char^{\lambda +\Lambda}$ is in a canonical bijection with the union over the two sided cells of $W_\zeta$
of the sets of pairs $\{ (\gamma, \psi)\}/{\mathcal G}_\uc$ where $\gamma\in {\mathcal G}_\uc$ and $\psi$ is
an irreducible representation of the centralizer of $\gamma$ in ${\mathcal G}_\uc$.

b) The category $(\Char^\hLambda_\ucs,\circ)$ is canonically equivalent to the category
$Sh^{\GG_\uc}(\GG_\uc)$ of conjugation equivariant sheaves of finite dimensional vector
spaces on $\GG_\uc$. This equivalence has the structure of a quasitensor equivalence.
If the cell $\uc$ is not exceptional,\footnote{Recall that there are no
exceptional cells unless $G$ has a factor of type $E_7$ or $E_8$. If $W_\zeta$
is of type $E_8$ there are two exceptional cells, and if it is of type $E_7$ there is one such cell.} it is moreover a tensor equivalence.
\end{cor}

\proofpt We start with part (b).

In \cite{BFO} we have shown, under the assumption (!), that for each 2-sided cell $\uc$ there exists a finite set $X_\uc$ equipped with a $\GG_\uc$ action and an equivalence between the
truncated convolution category ${\mathcal C}_\zeta^{\uc}$
\cite[\S 6]{BFO}  and the category
$Sh^{\GG_\uc}(X_\uc\times X_\uc)$ of $\GG_\uc$ equivariant sheaves of finite dimensional vector spaces on $X_\uc^2$. Furthermore, the category $Sh^{\GG_\uc}(X_\uc\times X_\uc)$
carries a natural monoidal structure given by convolution (which is a categorical counterpart of the matrix multiplication).
We let $\star$ denote the natural tensor structure on $Sh^{\GG_\uc}(X_\uc\times X_\uc)$
with associativity isomorphism {\em twisted} by the $\{\pm 1\}$ valued 3-cocycle appearing in
\cite[Proposition 4]{BFO}; the twisting is only nontrivial for exceptional cells.
The equivalence
${\mathcal C}_\zeta^{\uc} \cong Sh^{\GG_\uc}(X_\uc\times X_\uc)$ of \cite{BFO} is monoidal
with respect to the truncated convolution monoidal structure on the first category and $\star$ product on the second one.

 By Corollary~\ref{trunc}.b) the truncated convolution category ${\mathcal C}_\zeta^{\uc}$ is monoidally equivalent to $\Har^\hlambda_\ucs$.
Thus Theorem \ref{december} implies that $(\Char^\hlambda_\ucs,\circ)$
is tensor equivalent to the center of $(Sh^{\GG_\uc}(X_\uc\times X_\uc),\star)$.

By a result of \cite{O}, for a finite group $\Gamma$ and a finite set $X$ equipped with a $\Gamma$ action, the center
of the category $Sh^\Gamma(X\times X)$ (equipped with its natural monoidal structure) is canonically equivalent to the category
of conjugation equivariant sheaves $Sh^\Gamma(\Gamma)$ (in particular, it does not depend on $X$, this is a ``categorification" of the fact that centers of Morita equivalent associative
algebras are isomorphic; another application of this principle is in section
\ref{Mor_inv} above). Thus statement (b) is proven for non-exceptional cells.
For the monoidal category $Sh^\Gamma(X\times X)$ with associativity constraint
twisted by a 3-cocycle with values in $\{ \pm 1\}$, the center is quasitensor equivalent to the category $Sh^\Gamma(\Gamma)$, see e.g.
\cite{CGR}. This proves statement (b) in the remaining cases.

 Since irreducible objects in
 $Sh^{\Gamma_\uc}(\Gamma_\uc)$ are in bijection with pairs
 $(\gamma\in {\mathcal G}_\uc, \psi)$ as above, to prove (a) it remains to check that for
each irreducible
 character sheaf $\F$ there exists a unique cell $\uc$ such that
$\F\in \Char^\hLambda_{\leq \uc},\ \F\not \in \Char^\hLambda_{<\uc}$.

 Let $M$ be an irreducible submodule in $\Gamma^{\hlambda}(G,\F)$ and let
 $\uc$ be the cell of $M$. Then
 we have a nonzero map $\Ind_{U\otimes U}^{{\mathcal D}(G)}(M)\to \F$
which is surjective since $\F$ is irreducible.
 Thus $\F\in \Char^\hLambda_{\leq \uc}$. On the other hand $\F\not \in
\Char^\hLambda_{<\uc}$ since $\Har^\hlambda_{<c}\not \owns M\subset\Gamma^{\hlambda}(G,\F)$. \qed

\subsection{Semi-simplicity of the reduced cell category} 
\label{recL}

Recall the irreducible $U$-modules $L^\lambda_w$ introduced in~\ref{dual}. Let
$\fI_\uc=\bigcap_{w\in\uc}Ann(L^\lambda_w) $ be the intersection of primitive ideals attached to the cell
$\uc$.

Let $\Har_{\leq \uc,red}^\hlambda\subset \Har_{\leq \uc}^\hlambda$ be the full subcategory of bimodules
whose left and right annihilators contain $\fI_\uc$.

The following claim follows from \cite[Corollary 1.3.3]{Los}.

\begin{prop}\label{semisi} The image of the projection functor $\Har_{\leq \uc,red}^\hlambda\to
\Har_\uc^\hlambda$ consists of semi-simple objects. \end{prop}

\subsection{Proof of Theorem \ref{december}} \label{54} We construct the  functor in the opposite direction
$\fB_\uc: \Char^\hLambda_\ucs\to\CZ(\Har^\hlambda_\ucs)$.

Let $I_\uc$ be the unit object in $\Har^\hlambda_\ucs$; recall that it can be described as the sum of quotients
of $U(\g)$ by the primitive ideals belonging to the cell, or as the sum of irreducible objects
corresponding to Duflo (distinguished) involutions in the cell \cite{Ltc}.

 We set $\fB_\uc (\F) = \Gamma^\hlambda(G,\F) *
I_\uc\mod \Har^\hlambda_{<\uc}$.

By centrality of $\HC(\F)$ we have $\Gamma^\hlambda(G,\F) * I_\uc \cong I_\uc * \Gamma^\hlambda(G,\F)$, thus
this module is annihilated by $\fI_\uc$ both on the left and on the right. Thus by Proposition \ref{semisi}
the object $\fB_\uc(\F)$ is semi-simple.

The central structure of $\Gamma^\hlambda(G,\F)$ induces a central structure on $\fB_\uc(\F)$ via the
isomorphisms $$L*\fB_\uc(\F) \cong  L* \Gamma^\hlambda(G,\F)  \mod \Har^\hlambda_{<\uc} \cong
 \Gamma^\hlambda(G,\F) *L \mod \Har^\hlambda_{<\uc}  \cong
 \fB_\uc(\F) *L. $$

Here $*$ stands for tensor product in $\Har^\hlambda_\uc$ and $\Har^\hlambda$.

We proceed to check that
\begin{equation}\label{BAId}
\fB_\uc \circ \fA_\uc \cong\Id_{\CZ(\Har_\ucs^\hlambda)}.
\end{equation}

Recall the following general statement.

\begin{lem}\label{commz} (see \cite[Proposition 5.4]{ENO}) Let $\A$ be a semisimple rigid monoidal
category with a finite number of isomorphism classes of irreducible objects. Let $L_i$ be a set of
representatives for isomorphism classes of irreducible objects.
Let $F:\CZ(\A)\to \A$ be the forgetful functor, and $I:\CZ(\A)\to \A$ be its right adjoint.

Then the functor $F\circ I : \A \to \A$ is
naturally isomorphic to the functor $Av: M\mapsto \oplus_i L_i M L_i^*$. In particular

a) For $M\in \A$ the object $\oplus_i L_i M L_i^*$
 has a natural structure
of a central object.

b) The functor $Av: M\mapsto \oplus_i L_i M L_i^*$  is a commutator functor $\A\to \CZ(\A)$.

c) Moreover, $Av$ is a universal commutator functor from $\A$, i.e. for an additive category $\B$ the
category of functors $\CZ(\A)\to \B$ is equivalent to the category of commutator functors $\A\to \B$, where
the equivalence sends $\phi$ to $\phi\circ Av$. \end{lem}


It will be convenient to restate Lemma \ref{commz} in more functorial terms.

For a $\BC$-linear abelian category $\A$ 
 we can form a new abelian category $\A\otimes_\BC \A^{op}$.
The bilinear functor from $\A\times \A^{op}\to End(\A)$, $(M,N)\mapsto F_{M,N}:X\mapsto Hom(N,X)\otimes M$
induces a functor $\Phi: \A\otimes_\BC \A^{op}\to End(\A)$;
here $End(\A)$ denotes the category of  $\BC$-linear functors
 from $\A$ to $\A$. If $\A$ is semisimple with a finite number of irreducible objects,
 the functor $\Phi$ is easily seen to be an equivalence. 
Let $R_\A\in \A\otimes_\BC \A^{op}$ be defined by $R_\A=\Phi^{-1}(Id_{\A})$.

Furthermore, a monoidal structure on $\A$ induces a monoidal structure on
$\A\otimes_\BC \A^{op}$, and, if $\A$ is rigid, we
have a ``regular bimodule" action $\rho$ of $\A\otimes_\BC \A^{op}$ on $\A$,
$\rho(M\boxtimes N)(X)=M\otimes
X \otimes N^*$ (where we used rigidity to identify $\A$ and $\A^{op}$).

Also, there are two actions $\alpha_1$, $\alpha_4$ of the monoidal category $\A$ on $\A\otimes_\BC \A^{op}$
characterized by $\Phi(\alpha_1(M)(B))=L(M)\circ \Phi(B)$, $\Phi(\alpha_4(M)(B))=\Phi(B)\circ R(M)$ where
$L(M)$, $R(M)$ denote  left and right multiplication by $M$ in the monoidal category $\A$.
Likewise, there are two actions of $\A^{op}$ on $\A\otimes_\BC \A^{op}$ characterized by
$\Phi(\alpha_2(M)(B))= \Phi(B)\circ L(M^*)$, $\Phi(\alpha_3(M)(B))=R(M^*)\circ \Phi(B)$.

The tautological isomorphisms $L(M)\circ Id \cong Id\circ L(M)$ induce an isomorphism $\alpha_1(M)(R_\A)\cong
\alpha_4(M)(R_\A)$ for any $M\in \A$, which yields a central structure on the functor $\rho(R_\A):\A\to \A$.
Likewise, the isomorphism $R(M)\circ Id\cong Id\circ R(M)$ yields a commutator structure on $\rho(R_\A)$.

Lemma \ref{commz} can now be reformulated as follows

\begin{lem}\label{commz1} In notations and assumptions of Lemma \ref{commz}, we have an isomorphism
of functors $F\circ I\cong \rho(R_\A)$ compatible with the commutator and the central structures.

Also, $\rho(R_\A)$ is the universal commutator functor from $\A$. \end{lem}

We now return to the case $\A=\Har_\ucs^\hLambda$. In view of Lemma \ref{commz1} verification
of \eqref{BAId} reduces to construction of an isomorphism of endo-functors of $\Har_\ucs^\hLambda$
\begin{equation}\label{LML}
\fB_\uc \circ \fA_\uc\cong \rho(R_{\Har_\ucs^\hLambda})\end{equation}
compatible with the central and the commutator structures. (Recall that we abuse notation by denoting
a functor to the center of a category and its composition with the forgetful functor in the same way).

It is easy to see that $\A\otimes \A^{op}$ is identified with $\Har_{\uc\times \ucs}^{{\widehat{(\lambda,\lambda)
+\Lambda^2}}} $,
the semi-simple cell subquotient category in Harish-Chandra bimodules for the group $G\times G$.
The functor $\Phi$ is given by $\Phi(B)(M)=B\otimes _{U_3\otimes U_4} M \mod \Har^\hLambda_{<\uc}$,
while the action $\rho$ is given by $\rho(B)(M)=B\otimes _{U_2\otimes U_3}M$ in the self-explanatory notation.

Thus \eqref{LML} would follow if we show that $\fB_\uc \circ \fA_\uc:
M\mapsto R\otimes _{U_2\otimes U_3}M \mod
\Har^\hLambda_{<\uc}$
for some quadro-module $R_\uc$ such that $R_\uc\otimes _{U_3\otimes U_4} M\cong M$ for any $M\in
\Har^\hLambda_{\uc}$; moreover, the central and the commutator structure on $\fB_\uc \circ \fA_\uc$
should be compatible with the one coming from the last isomorphism.

Recall that $ \Gamma^\hlambda(G,-) \circ \Ind:  
M\mapsto \CalD(G)\otimes _{U_2\otimes U_3} M$,
where $ \CalD$ is equipped with four pairwise commuting $\g$-actions coming, respectively,
from left multiplication by left invariant vector fields, right multiplication by left invariant vector fields,
right multiplication by right invariant vector fields and left multiplication by right invariant vector fields.



 Notice that the formula $\CalD(G) = Ind_{U_{3,4}}^{U_3\otimes U_4}({\mathcal O}(G))$ 
where $U_{3,4}\subset U_3\otimes U_4$ is the image of the diagonal homomorphism
and $U$ acts on ${\mathcal O}(G)$ by differentiation with respect to right invariant vector fields,
yields a functorial isomorphism
\begin{equation}\label{DId} \CalD(G) \otimes _{U_3\otimes U_4}M\cong M \end{equation}
for a Harish-Chandra bimodule $M$.

Thus the quotient $\CalD(G)_{\leq \uc} := \CalD(G)
\otimes_{U^{\otimes 4}} (U/\fI_\uc)^{\otimes 4}$
satisfies $\CalD(G)_{\leq \uc}\otimes_{U_3\otimes U_4} M\cong M$ for any $M$ annihilated by $\fI_\uc$.
It follows that the object $R_\uc:=\CalD(G)_{\leq \uc}\mod \Har_{<\uc\times \uc}^{{\widehat{(\lambda,\lambda)
+\Lambda^2}}} $ satisfies the characterizing property of $R_{\Har_\uc^\hLambda}$. It is not hard to check
that the centrality and the commutator structure on $\fB_\uc\circ \fA_\uc$ are compatible with those
coming from the isomorphism of functors $R_\uc\otimes_{U_3\otimes U_4} - \cong Id_{\Har_\uc^\hLambda}$.
This proves \eqref{LML}, and hence \eqref{BAId}.

 We also claim that the functor $\fA_\uc$ is semi-simple,
while $\fB_\uc$ is 
conservative (i.e. it kills no objects). Semisimplicity of $\fA_\uc$ follows from Decomposition
Theorem, since $\CH$ is a composition of smooth pull-back and proper push-forward (see \ref{CSdef}).
Conservativity of
$\fB_\uc$ comes from the following observation:
 the definition of cells in character sheaves (see \ref{cCS}) easily implies that for
 $M\in \Char_{\leq \uc}^\hLambda$, $M\not \in \Char_{< \uc}^\hLambda$ there exists a surjection
in $\Har^\hlambda/\Har^\hlambda_{<\uc}$:
 $\Gamma^\hlambda(M)\to L$ for some irreducible Harish-Chandra bimodule $L$ belonging to $\uc$. Since
 $L$ is annihilated (say, on the right) by $\fI_\uc$, we see that $L\otimes_{U\otimes U} I_\uc = L$, so
$\Gamma^{\hlambda}(M)\otimes_{U\otimes U}I_\uc$ surjects to $L=L\otimes_{U\otimes U} I_\uc $ in
$\Har^\hlambda/\Har^\hlambda_{<\uc}$, in particular $\Gamma^{\hlambda}(M)\otimes_{U\otimes U}I_\uc
\ne 0 \mod \Har^\hlambda_{<\uc}$.

These properties of $\fA_\uc$, $\fB_\uc$ together with \eqref{BAId} imply that $\fA_\uc$ sends
irreducible objects to irreducible ones and induces an injective map between the sets of isomorphism classes
of semi-simple objects. To show that this map is surjective assume, on the contrary, that an irreducible
object $L\in \Char_\uc^\hLambda$ does not lie in the image of $\fA_\uc$. 
Then $Hom_{ \Char_\uc^\hLambda}(\fA_\uc(Z),L)=0$ for any
$Z\in \CZ(\Har_\ucs^\hlambda)$. 
This is equivalent to saying that $Hom_{ \Char_\uc^\hLambda}(
Ind (M), L)=0$ for any $M\in \Har_\ucs^\hlambda$
(where we use notations introduced before Proposition \ref{exact}).
Applying adjointness pointed out before Proposition \ref{exact} we get
that $Hom_{\Har_\uc^\hlambda}(M, \Gamma^\hlambda(G,L) )=0$ for any $M\in \Har_\ucs^\hlambda$.
This implies that $\Gamma^\hlambda(G,L)\in \Har_\uc^\hlambda$ vanishes, which 
 contradicts conservativity of $\fB_\uc$.

Thus we have proven Theorem \ref{december}(a). It remains to equip the equivalence $\fA_\uc$
(equivalently, $\fB_\uc$) with a tensor
structure. The functor $\HC:\ D\Char^\hLambda\to\CZ(D\Har^\hlambda)$ is monoidal. The object $I_\uc$
satisfies $$I_\uc\Lotimes_U I_\uc\cong I_\uc \oplus R \mod D\Har^\hlambda_{<\uc}$$ where $R\in D^{<0}
(\Har^\hLambda)$, this is so since the convolution of semi-simple objects is a semi-simple complex and
$I_\uc$ is a unit object with respect to truncated convolution on the cell subquotient category.

So for $\F,\G\in \Char^\hLambda_\ucs$ we have: \begin{multline*} \fB_\uc(\F\circ \G)\cong\HC( H^0(\F*\G))\Lotimes _U
I_\uc \mod {\Har^\hlambda_{<\uc}} \cong H^0(\F*\G*I_\uc) \mod \Har^\hlambda_{<\uc} \cong \\
H^0(H^0(\F*I_\uc)\Lotimes_U H^0(\G*I_\uc)) \mod \Har^\hlambda_{<\uc}\cong \fB_\uc(\F)\otimes_U\fB_\uc(\G).
\end{multline*}

This endows $\fB_\uc$ with a monoidal structure thereby completing the proof of Theorem \ref{december}.

\section{Specialization of character sheaves}\label{6} In this section we sketch a proof
 that the exact functor of Corollary \ref{exa} can
 be described as Verdier
specialization in the De Concini -- Procesi compactification (in particular, it commutes with Verdier
duality).

Let $G^{ad}$ denote the quotient of $G$ by its center and $T^{ad}$ its maximal torus. Set
$\CX^{ad}:=(G^{ad}/N^+\times G^{ad}/N^-)/T$, let $\pi:G\to G^{ad}$, $\pi^\CX:\CX\to \CX^{ad}$ denote the
 projections.

Let $\ol G^{ad}$ denote De Concini -- Procesi compactification of $G^{ad}$ \cite{DCP}, see also exposition in
\cite{EJ}; let $j:G^{ad}\imbed \ol G^{ad}$ be the embedding. The complement $Z = \ol G^{ad}\setminus G^{ad}$
is a divisor with normal crossing. The closed stratum (intersection of the components) in
 $Z$ is identified with
$\CB^2$ and the complement in the total space of the normal bundle $N_{\CB^2}(\ol G^{ad})$ to the normal
cone $N_{\CB^2}(Z)$ is identified with $\CX^{ad}$.

We have  the Verdier specialization functor
 $ Sp:\ \CalD{\mathrm -mod}_{rs}(\ol G^{ad})\to\CalD
{\mathrm -mod}(N_{\CB^2}(\ol G^{ad})$. It induces a functor $ Sp^o:\ \CalD{\mathrm -mod}_{rs}
(G^{ad})\to\CalD{\mathrm -mod}(\CX^{ad})$. Here the subscript $"rs"$ refers to the category of holonomic
$\CalD$-modules with regular singularities.

\begin{thm}\label{glsp} Let $M$ be a holonomic $D$-module with regular singularities
 on $G^{ad}$ such that the action of the center $Z(U(\g))$
coming from the action of $G^{ad}$ on itself by left translations is locally finite. Then the $D$-module
$Sp(j_*M)$ is $T^{ad}$-monodromic and for a dominant weight $\lambda$ we have a canonical isomorphism:
$$\Gamma^{\hat \lambda}(Sp^o(M))\cong \Gamma^{\hat \lambda}(M).$$ \end{thm}

\begin{cor}\label{spec}
 For a character sheaf $A$ on $G$ we have a
 canonical isomorphism
 $$ Sp^o(\pi_*A)\cong \pi^\CX_*\,  \bI_{w_0}\circ\HC(A).$$
\end{cor}

\proof of Corollary. The isomorphism is immediate from Theorem \ref{glsp} and Proposition \ref{HotKas}.
\qed

\proof of the Theorem. Let $\CI$ be the ideal sheaf of $\CB^2\subset \ol G^{ad}$. Let $\CalD_\CI\subset
\CalD$ be the subsheaf of differential operators preserving each power $\CI^n$, $n\in \Zet_{>0}$. The
associated graded sheaf $\ol \CalD_\CI$ acts faithfully on  $\oplus \CI^n/\CI^{n+1}$.

The set of components of $Z$ is in a canonical bijection with the set $I$ of vertices of the Dynkin diagram
of $G$.
 For $i\in I$ let $Z_i$ be the corresponding component and $\alpha_i$
the corresponding simple root, $\check{\alpha}_i$ be the simple coroot. Since $Z_i$ are smooth divisors and
$\CB^2$ is their transversal intersection, the conormal bundle $\CI/\CI^2$ is decomposed as a direct sum of
line subbundles indexed by $I$. We have a canonical embedding $e:\ft\to \Gamma (\ol \CalD_\CI)$ such that
$x\in \ft$ acts on the $i$-th line subbundle by $\alpha_i(x)$; we will use the same notation for the induced
map $\on{Sym(\ft)}\to \Gamma(\ol\CalD_\CI)$. Let $E=e(\check{\rho})$ be the Euler vector field.

It is known that a nonzero
 function $f\in V_\lambda\otimes V_\lambda^*\subset \CO(G^{ad})$
has a pole of order $\lambda(\check{\alpha}_i)$ along the divisor $Z_i$. Hence we have
\begin{equation}\label{ehc} e|_{\on{Sym}(\ft)^W}=\ol{hc}, \end{equation} where
$\ol{hc}:\on{Sym}(\ft)^W=Z(U(\g))\to \Gamma(\ol{\CalD_\CI})$ comes from the (left) action of $G$ on $\ol
G^{ad}$.

Recall that the specialization of a $\CalD$-module is defined as the associated graded of the
Malgrange-Kashiwara $V$-filtration \cite{K}. By the definition, each term of the $V$-filtration
 is invariant with respect to $ \CalD_\CI$
and the action of the Euler vector field $E$ on $gr^d(V^\bullet)$ is finite with all eigenvalues in $S+d$
where $S\subset \Ce$ is a fixed set of representatives for the cosets $ \Ce/\Zet$.

Thus \eqref{ehc} shows that for $M$ as in the Theorem the action of $\ft$ on $gr(\bV^\bullet(j_*(M))$ is
locally finite. Hence $Sp(M)$ is a $T^{ad}$-monodromic (weakly equivariant) $\CalD$-module.

Furthermore, it follows that the filtration on $\Gamma(M)=\Gamma(j_*(M))$ induced by the $V$-filtration is
given by $$\bV^i(\Gamma(M))=\oplusl_{\lambda | \lambda(\check{\rho})\in S+d, d< -i}
\Gamma^{\hat{\lambda}}(M).$$

Thus we see that
 $$\Gamma^{\hat{\lambda}} (Sp(j_*(M))=
\Gamma^{\hat{\lambda}}(gr(\bV^\bullet(M)))\subset gr (\bV^\bullet \Gamma (j_*M))^{\hat\lambda}=\Gamma^{\hat
\lambda}(M),$$ i.e. we get a canonical injective map $i_\lambda$ from the right hand side to the left hand
side of the desired isomorphism.

To check that this map induces an isomorphism for dominant $\lambda$ consider a surjective map of
quasicoherent sheaves on $\ol G$: $s_\lambda:(j_*M)^{\leq \lambda}\to gr \bV(j_*M)^{\hat \lambda}$. Here $
gr \bV(j_*M)^{\hat \lambda}$ is the sheaf supported on $\CB^2$ which is the maximal direct summand in
$gr^i\bV(j_*M)$ where $\ft$ acts by generalized character $\lambda$; here $i$ is such that
$\lambda(\check{\rho})\in i+S$. The sheaf $(j_*M)^{\leq \lambda}$ is the preimage of $ gr \bV(j_*M)^{\hat
\lambda}$ under the map $\bV ^i(j_*M)\to gr^i\bV(j_*M)$.

For an element $\nu=\sum n_i\alpha_i$ in the root lattice let $Z^\nu$ denote the divisor $\sum n_iZ_i$ on
$\ol G$. It is known that $O(Z^\nu)|_{\CB^2}\cong \CO(\nu,\nu)$. The line bundle $\CO(Z^\nu)$ is ample
provided that the weight $\nu$ is dominant. It follows from the definitions that for $\mu=\lambda+\nu$ with
$\nu$ in the root lattice, the map $s_\mu$ is obtained from $s_\lambda$ by twisting with $\CO(Z^\nu)$.

Take $\sigma\in \Gamma^{\hat \lambda}(Sp(j_*M))=\Gamma(gr \bV(j_*M)^{\hat \lambda})$. For $\nu$ positive
enough the sections of $\CO_{\CB^2}(\nu,\nu)\subset gr \bV(j_*M)^{\hat {\lambda+\nu}}$ (where the embedding
comes from $\sigma$) lie in the image of $\Gamma(s_{\lambda+\nu})$. However, the image of
$\Gamma(s_\lambda)$ coincides with the image of $i_\lambda$. One can show that for  dominant weights
$\lambda,\, \mu$ with $\lambda-\mu$ in the root lattice, the $\g\oplus \g$ module  $\Gamma^{\hat \mu}(M)$
is obtained
 $\Gamma^{\hat \lambda}(M)$ by translation functor.
It follows that the maps $i_\lambda$, for various $\lambda$, come from a fixed map of $T^{ad}$-monodromic
$\CalD$-modules on $\CX$. Thus Localization Theorem implies that $\sigma$ is in the image of $i_\lambda$.

We have checked that $ \Gamma^{\hat \lambda}(M)\cong \Gamma^{\hat \lambda}(Sp (j_*M))$. However, it is not
hard to deduce from Kashiwara Lemma that for a $T^{ad}$-monodromic $\CalD$-module $N$ on $N_{\CB^2}(\ol G)$
and a dominant weight $\lambda$ we have $\Gamma^{\hat \lambda}(N)=\Gamma^{\hat \lambda}(N|_{\CX^{ad}})$.
\qed

\footnotesize{ {\bf R.B.}: Department of Mathematics, Massachusetts Institute of Technology,\\
Cambridge, MA 02139, USA;\\ 
{\tt bezrukav@math.mit.edu}}

\footnotesize{ {\bf M.F.}: IMU, IITP, and National Research
University Higher School of Economics,\\ Department of
Mathematics, 20 Myasnitskaya st, Moscow 101000 Russia;\\ 
{\tt fnklberg@gmail.com}}

\footnotesize{ {\bf V.O.}: $\,$  Department of Mathematics, 1222 University of Oregon, Eugene OR 97403-1222,
USA\\ 
{\tt vostrik@math.uoregon.edu}}

\end{document}